\documentstyle{article}
\fussy
\newtheorem{theorem}{Theorem}[section]
\newtheorem{lm}{Lemma}[section]
\newtheorem{prop}{Proposition}[section]

\def\k{{K\"{a}hler }}

\def\cy{{Calabi-Yau }}
\def\l{{Lagrangian }}

\input epsf
\begin{document}
\hbadness=10000
\title{\bf Lagrangian Torus Fibrations and Mirror Symmetry of Calabi-Yau Manifolds}
\author{Wei-Dong Ruan\\
\it Department of mathematics\\
\it University of Illinois at Chicago\\
\it Chicago, IL 60607\\
\it Email: ruan@math.uic.edu}
\date{}
\footnotetext{Partially supported by NSF Grant DMS-9703870 and AMS Centennial Fellowship.}
\maketitle
%\tableofcontents
\section{Introduction}
A very beautiful part of Kodaira's complex surface theory is the theory of elliptic surfaces. Elliptic surfaces lie somewhere between ``positive'' (rational) surfaces and ``negative'' (general type) surfaces, where many wonderful things happen. Kodaira's elliptic surface theory, especially his classification of singular fibres, have lots of deep relations with other parts of mathematics, including the theory of elliptic curves and singularity theory. In recent years, a rather surprising higher dimensional generalization of such classical and fundamental nature to Calabi-Yau manifolds was discovered in the context of mirror symmetry!\\

Mirror Symmetry conjecture originated from physicists' work in conformal field theory and string theory. It proposes that for a Calabi-Yau 3-fold $X$ there exists a Calabi-Yau 3-fold $Y$ as its mirror. The quantum geometry of $X$ and $Y$ are closely related. In particular one can compute the number of rational curves in $X$ by solving the Picard-Fuchs equation coming from variation of Hodge structure of $Y$ \cite{Can}. Before mirror symmetry,  Calabi-Yau manifolds were rather abstract and mysterious in mathematics. The study of mirror symmetry, especially the exciting development of counting of rational curves gave us a lot of insight of quantum geometry of Calabi-Yau manifolds or more precisely the geometry of the moduli spaces of Calabi-Yau manifolds. However the understanding of classical geometry and topology of general Calabi-Yau manifolds and how mirror symmetry is reflected in the classical level was still lacking. Strominger-Yau-Zaslow's mirror conjecture, which is very much in the spirit of Kodaira elliptic surface theory, changed all these.\\

In 1996 Strominger, Yau and Zaslow (\cite{SYZ}) proposed a geometric construction of mirror manifold via special Lagrangian torus fibration. According to their program (we will call it SYZ construction), a Calabi-Yau 3-fold should admit a special Lagrangian torus fibration. The mirror manifold can be obtained by dualizing the fibres. Or equivalently, the mirror manifold of $X$ is the moduli space of special Lagrangian 3-torus in $X$ with a flat $U(1)$ connection.\\

Notice that despite its physical root, the statement of SYZ conjecture is purely mathematical and rather classical in nature. It makes the structure of Calabi-Yau manifolds rather explicit and mirror symmetry a rather concrete duality. In the K3 surface case, K3 with SYZ fibration is equivalent to elliptic K3 under the hyper\k twist.\\
 
However, in spite of the great promise of SYZ conjecture, our understanding on special Lagrangian submanifolds and therefore the SYZ fibration is very limited. The known examples are mostly explicit local examples or examples coming from 2-dimensional case. There are very few examples of special Lagrangian submanifold or special Lagrangian fibration for dimension higher than two. M. Gross, P.M.H. Wilson and N. Hitchin (\cite{Gross1}\cite{Gross2}\cite{GW}\cite{H}) did some important work in this area in recent years. These works mainly concern local geometric structure of the special Lagrangian fibration and cases that can be reduced to 2-dimensional situation. On the other extreme, in \cite{Z}, Zharkov constructed some non-Lagrangian torus fibration of Calabi-Yau hypersurface in toric variety. Despite all these efforts, SYZ construction still remains a beautiful dream to us.\\
 
Given the general lack of knowledge for special Lagrangian, our approach to SYZ conjecture is to relax the special Lagrangian condition and consider Lagrangian fibration, which we feel like to be a good compromise and is interesting in its own. Special Lagrangians are very rigid, while Lagrangian submanifolds are more flexible and can be modified locally by Hamiltonian deformation. However for many application to mirror symmetry, especially those concerning (symplectic) topological structure of fibration, Lagrangian fibration will provide quite sufficient information. In our work \cite{lag1, lag2, lag3, tor}, we mainly concern Lagrangian torus fibrations of Calabi-Yau hypersurfaces in toric variety, namely the symplectic topological aspect of SYZ mirror construction.\\

In light of the recent important explicit local examples of generic special Lagrangian fibrations constructed by Dominic Joyce \cite{Joyce}, the SYZ fibration map of Calabi-Yau manifolds are likely to be only piecewise $C^\infty$ (Lipschitz) instead of $C^\infty$, and the singular locus instead of being a codimension two graph in $S^3$ is likely to be of codimension one as a fattening of a graph. This is very much in line with the codimension one singular locus pictures that naturally come out of our gradient flow construction in \cite{lag1}. As pointed out by Joyce, the actual special Lagrangian fibration probably can be constructed by perturbing our Lagrangian fibration with codimension one singular locus. SYZ duality of the fibres are likely to be precise duality only at the large complex (radius) limit. In my point of view, this does not necessarily destroy the beauty of SYZ construction, rather makes it much richer and more intriguing. For the sake of computation, understanding the geometric and topological structure of Calabi-Yau manifolds, especially, constructing the mirror Calabi-Yau symplectic topologically, it is still rather convenient to have a Lagrangian fibration with codimension two graph singular locus.\\

In \cite{tor}, we were able to construct Lagrangian torus fibrations with codimension two graph singular locus of generic Calabi-Yau hypersurfaces in toric varieties correponding to reflexive polyhedra in complete generality. With these detailed understanding of Lagrangian torus fibrations of generic Calabi-Yau hypersurfaces in toric varieties, we were able to prove the symplectic topological version of SYZ mirror conjecture for Calabi-Yau hypersurfaces in toric varieties, although the origional SYZ conjecture for special Lagrangian fibrations probably need some modification. More precisely we have:\\
\begin{theorem}
For generic Calabi-Yau hypersurface $X$ in the toric variety corresponding to a reflexive polyhedron $\Delta$ and its mirror Calabi-Yau hypersurface $Y$ in the toric variety corresponding to the dual reflexive polyhedron $\Delta^\vee$ near their corresponding large complex limit and large radius limit, there exist corresponding Lagrangian torus fibrations\\
\[
\begin{array}{ccccccc}
X_{\phi(b)}&\hookrightarrow& X& \ \ &Y_b&\hookrightarrow& Y\\
&&\downarrow& \ &&&\downarrow\\
&& \partial \Delta_v& \ &&& \partial \Delta_w\\
\end{array}
\]\\
with singular locus $\Gamma \subset \partial \Delta_v$ and $\Gamma' \subset \partial \Delta_w$, where $\phi:\partial \Delta_w \rightarrow \partial \Delta_v$ is a natural homeomorphism that satisfies $\phi(\Gamma')=\Gamma$. For $b\in \partial \Delta_w \backslash \Gamma'$, the corresponding fibres $X_{\phi(b)}$ and $Y_b$ are naturally dual to each other.\\
\end{theorem}
Our work essentially indicates that the Batyrev-Borisov mirror construction \cite{Bat,Bor}, which was proposed purely from toric geometry stand point, can also be understood and justified by the SYZ mirror construction. This should give us greater confidence on SYZ mirror conjecture for general Calabi-Yau manifolds.\\

In an earlier physics work of Leung and Vafa \cite{lv}, they discussed a heuristic derivation of Batyrev's mirror construction via T-duality. Although we do not fully understand the physics argument presented in \cite{lv}, we later find that some ideas in \cite{lv} seem to be closely related to our approach toward SYZ mirror symmetry in the Batyrev case and may be helpful to better understand the relation of our mathematical construction and its underlying physics.\\
  
Around the same time as our paper \cite{lag3}, in which Lagrangian torus fibrations were constructed for generic quintic Calabi-Yau and corresponding symplectic SYZ was proved, M. Gross independently (using completely different method) constructed certain {\it non-Lagrangian} torus fibrations for quintic Calabi-Yau in \cite{Gross3}, which exhibit similar topological structure as our Lagrangian torus fibration in \cite{lag3}. His idea, very much in line with part of our previous work \cite{lag1}, was to use the monodromy information to guess the structure of the fibration topologically. (Here he got the expected monodromy information from the resolution of orbifold singularities in the mirror.) He then constructed the topological manifold accordingly and used a theorem of C.T.C. Wall to prove the topological manifold constructed is homeomorphic to quintic Calabi-Yau, therefore constructing a topological torus fibration on the quintic Calabi-Yau. On the other hand, our completely different approach in \cite{lag3} constructed the Lagrangian torus fibration on quintic Calabi-Yau directly (without going to the mirror) using gradient flow based on the graphlike behavior of ameoba/string diagram discussed in \cite{N,M}. We note that our approach (in \cite{lag3}) also provides an alternative way (from \cite{Gross3}) to construct {\it topological} torus fibrations on quintic Calabi-Yau directly. The technical results required for such alternative topological construction will be much easier than the technical results required for the construction of Lagrangian torus fibration discussed in \cite{lag2,N}.\\

The purpose of this paper is to give a coherent overview of our work on Lagrangian torus fibration and symplectic SYZ mirror correspondence in various cases, starting from the most famous case of quintic Calabi-Yau and its mirror. Through our discussion, we will pay special attention to exploring the relation between various changes of singular locus graphs of the fibrations and various geometric changes of the corresponding Calabi-Yau manifolds.\\

In the elliptic surface case, the singular locus is just a finite set of isolated points. (For a generic elliptic K3, there are $24$ points.) Generic type of singular fibre is the nodel ${\bf CP^1}$. Several singular locus points can come together in non-generic elliptic surfaces to form more complex singular fibres. Analogously for the case of Lagrangian torus fibrations of Calabi-Yau 3-folds over $S^3$, the singular locus is genericly a graph $\Gamma$ with only 3-valent vertex points $\Gamma^0$, which are seperated into two sets --- the positive (negative) vertex points $\Gamma^0_+$ ($\Gamma^0_-$)--- depending on whether Euler number of the singular fibre over it is $+1$ or $-1$ as described in \cite{tor}. The singular locus graph and singular fibres change from region to region in the moduli space by passing through some walls of non-generic Calabi-Yau's. The most generic of non-generic graphs in 3-space is a graph with only 3-valent vertex points except one 4-valent vertex. From the point of view of deforming graphs in 3-space, there are the following three natural ways to degenerate generic graphs into such non-generic graph as indicated in the following picture.\\
\begin{center}
\begin{picture}(200,170)(-150,30)
\thicklines
\put(-50,121){\line(0,1){24}}
\put(-50,121){\line(0,-1){24}}
\put(-50,121){\line(-1,0){30}}
\put(-50,121){\line(1,0){30}}
\put(-50,121){\circle*{4}}
\thinlines
\put(-100,135){\vector(-2,1){20}}
\put(-100,105){\vector(-2,-1){20}}
\put(-115,50){\thicklines
\put(-59,130){\line(0,1){24}}
\put(-49,120){\line(0,-1){24}}
\put(-59,130){\line(-1,0){30}}
\put(-49,120){\line(1,0){30}}
\put(-59,130){\line(1,-1){10}}
\put(-49,120){\circle*{4}}
\put(-59,130){\circle*{4}}}
\put(-115,-50){\thicklines
\put(-49,122){\line(0,1){24}}
\put(-59,112){\line(0,-1){24}}
\put(-59,112){\line(-1,0){30}}
\put(-49,122){\line(1,0){30}}
\put(-59,112){\line(1,1){10}}
\put(-49,122){\circle*{4}}
\put(-59,112){\circle*{4}}}
\thinlines
\put(-8,121){\vector(1,0){30}}
\put(125,-5){
\put(-57,128){\circle*{1}}
\put(-55,126){\circle*{1}}
\put(-53,124){\circle*{1}}
\put(-51,122){\circle*{1}}
\thicklines
\put(-49,120){\line(0,1){8}}
\put(-49,132){\line(0,1){20}}
\put(-49,120){\line(0,-1){20}}
\put(-59,130){\line(-1,0){30}}
\put(-59,130){\line(1,0){42}}}
\end{picture}
\end{center}
\begin{center}
\stepcounter{figure}
Figure \thefigure: Generic degeneration of graph
\vspace{10pt}
\end{center}

Notice that the graphs in the above picture are all located in parallel planes. (In our case planes parallel to the paper.) We will see the reason for this constrain later. There are two ways to collide two 3-valent points into one 4-valent point (as indicated on the left) and there is one way for two smooth lines to meet and form the 4-valent vertex. It turns out these graph operations canonically correspond to some natural geometric or topological change of Calabi-Yau manifolds. In the ground breaking work of Aspinwall, Greene and Morrison \cite{AGM2} on topological change in mirror symmetry, \k moduli of different birational equivalent Calabi-Yau as well as other non-Calabi-Yau regions are unified to form the \k moduli that can be identified with the complex moduli of the mirror Calabi-Yau via the monomial-divisor mirror map (\cite{AGM}). Different regions in the \k moduli corresponding to different birational Calabi-Yau models are related by flops. The transition between the two graphs on the left of Figure 1 in the \k side exactly corresponds to such a flop. Another way to extend the moduli space of Calabi-Yau is through a conifold transition (or black hole condensation in physics literature) that is related to the so-called Reid's fantasy \cite{reid}, which conjectures that all the moduli spaces of Calabi-Yau are connected through such transition. As suggested in \cite{black1, black2, black3, Mor}, the conifold transition can be used to find mirrors beyond toric cases. The transitions between graphs on the left and the graph on the right of Figure 1 exactly corresponds to a conifold transition. (Gross also observed such transition in \cite{Gross4}.) These correspondences further indicate the crucial importance of Lagrnagian torus fibration in understanding mirror symmetry. In Section 5 and 6, we will indicate these rather intriging relations through examples. We will also describe the corresponding changes of singular fibres.\\\\

\section{Quintic case}
Our idea of construction is a very natural one. We try to use gradient flow to get Lagrangian torus fibration from a known Lagrangian torus fibration at the "Large Complex Limit". This method will in principle be able to produce \l torus fibration in general \cy hypersurfaces or complete intersections in toric variety.\\

To illustrate our idea, it is helpful to explore the historically most famous case of quintic Calabi-Yau threefolds in ${\bf CP^4}$ in detail. Most of the essential features of the general case already show up here. Let\\
\[
z^m=\prod_{k=1}^5z_k^{m_k},\ |m|=\sum_{k=1}^5 m_k,\ {\rm for}\ m=(m_1,m_2,m_3,m_4,m_5)\in {\bf Z}^5_{\geq 0}.
\]
Then a general quintic can be denoted as\\
\[
p(z) = \sum_{|m|=5}a_mz^m.
\]\\
Let $m_0=(1,1,1,1,1)$, and denote $a_{m_0}=\psi$. Consider the quintic Calabi-Yau familly $\{X_\psi\}$ in ${\bf P^4}$ defined as\\
\[
p_\psi(z) = p_a(z) + \psi \prod_{k=1}^5 z_k = \sum_{m\not=m_0,|m|=5}a_mz^m + \psi \prod_{k=1}^5 z_k=0.
\]\\
When $\psi$ approaches $\infty$, the familly approach its ``Large Complex Limit'' $X_{\infty}$ defined by\\
\[
p_{\infty}=\prod_{k=1}^5 z_k=0.
\]\\
$X_{\infty}$ is a union of five ${\bf CP^3}$'s. There is a natural degenerate  $T^3$ fibration structure for $X_{\infty}$ given by the natural moment map $F: \ {\bf CP^4}\longrightarrow \Delta$.  $\Delta={\rm Image}(F)$ is a 4-simplex. $X_{\infty}$ is naturally fibered over $\partial \Delta$ via this map $F$ with general fibre being $T^3$. This is in a sense precisely the SYZ  special Lagrangian $T^3$ fibration for $X_{\infty}$.\\

Consider the meromorphic function\\
\[
s= \frac{p_\infty(z)}{p_a(z)} = f+ih
\]\\
defined on ${\bf CP^4}$. Let $\nabla f$ denote the gradient vector field of real function $f=Re(s)$ with respect to the \k metric $g$.\\

Notice that\\
\[
\nabla f = H_h,
\]\\
where $H_h$ is the Hamiltonian vector field generated by $h=Im(s)$. This implies the following\\
\begin{lm}
The gradient flow of $f$ leaves the set $\{Im(s)=0\}$ invariant and deforms \l submanifolds in $X_{\infty}$ to \l submanifolds in $X_{\psi}$.\\
\end{lm}

With this lemma in mind, the construction of \l torus fibration of $X_{\psi}$ for $\psi$ large is immediate. Deforming the canonical Lagrangian torus fibration of $X_{\infty}$ over $\partial \Delta$ along the gradient flow of $f$ will naturally induce a \l torus fibration of $X_{\psi}$ over $\partial \Delta$ for $\psi$ large and real.\\

The key advantage of the gradient flow method is that once it is run, the Lagrangian fibrations are automatically produced almost effortlessly. There are no adhoc manual construction involved to this point. Of course, for such Lagrangian fibration to be of any use, it is necessary to understand the detailed structure such as the singular fibres, singular set and the singular locus, etc. For this purpose, it is necessary to understand the details of dynamics of the gradient flow and how the fibrations evolve under the flow. There are a lot of rather delicate technical issues to be addressed here. First of all, the gradient flow in our situation is rather non-conventional. The critical points of the function are usually highly degenerate and often non-isolated. Worst of all the function is not even smooth (it has infinities along some subvarieties). In \cite{lag2} we discussed the local behavior of our gradient flow near critical points and infinities of $f$, and also the dependence on the metric to make sure that they behave the way we want.\\

Secondly, our gradient flow method naturally produces piecewise smooth (Lipschitz) Lagrangian fibrations. The singular locus is of codimension one. This is quite to the contrary of the conventional wisdom, where people expect the special Lagragian fibration to be smooth and the singular locus to be of codimension two. Nerverthless, in light of the recent examples of Joyce \cite{Joyce}, Lagrangian fibrations with codimension one singular locus might reflect the structure of the actual special Lagrangian fibration after all. This further indicates that it might be a good idea to consider Lagrangian fibrations, which for most purposes would be as good as, or even more convenient to use than the actural special Lagrangian fibrations. In \cite{lag2}, we squeeze the codimension one singular locus by symplectic geometry technique to get Lagrangian torus fibration with codimension two graph singular locus. For computation purposes and to construct mirror manifolds symplectic topologically, such models clearly are rather desirable. In particular, the SYZ duality can be made more precise in the symplectic category using such models. For details, please refer to \cite{lag2}.\\
 
With all the technical points taken care of, it is not hard to observe that the singular set of our fibration $F_\psi: X_\psi \rightarrow \partial \Delta$ is the curve\\
\[
C = X_\psi \cap {\rm Sing}(X_\infty).
\]\\
The singular locus\\
\[
\tilde{\Gamma} = F_\psi(C)
\]\\
is located in the 2-skeleton of $\Delta$, which is a union of 2-simplices. More precisely, since $C$ is reducible, and each irreducible component is in a ${\bf CP^2}$, according to \cite{lag2}, the study of singular locus can be isolated to each ${\bf CP^2}$ and reduced to the following problem.\\\\
{\bf Problem:} Let $F:{\bf CP^2} \rightarrow \Delta$ be the standard moment map with respect to the Fubini-Study metric. For what kind of quintic curve $C$ in ${\bf CP^2}$, $\tilde{\Gamma} = F(C)$ is a fattening of some graph $\Gamma$?\\

This problem was discussed in \cite{N} for curves $C$ of arbitrary degrees in ${\bf CP^2}$, more generally for curves in arbitrary 2-dimensional toric varieties. There is a notion of {\bf near the large complex limit}, which corresponds to the coefficients of the polynomial defining the curve satisfying some convexity condition with respect to the Newton polygon of the polynomial. For curves in ${\bf CP^2}$, we have\\
\begin{theorem}
When a degree $d$ curve $C_{p_d}$ defined by polynomial $p_d(z)$ is near the large complex limit, $F(C_{p_d})$ will have exactly $g = \displaystyle\frac{(d-1)(d-2)}{2}$ holes and $d$ external points in each edge of $\Delta$.
\end{theorem}

Different convexity condition $Z$ of the coefficients of $p_d(z)$ corresponds to different large complex limit chambers, which determine different graphs $\Gamma_Z\subset \Delta$, of which $\tilde{\Gamma} = F(C_{p_d})$ is a fattening. More precisely, we have\\
\begin{theorem}
For $t\in {\bf R}_+$ small, $F_{t^w}(C_{p_d})$ is a fattening of $\Gamma_Z$.\\
\end{theorem}
Here $t^w$ is the coefficient set of $p_d(z)$. $w$ satisfies the convexity condition $Z$. The theorem is roughly saying that when $p_d(z)$ is approaching the large complex limit corresponding to $Z$, $\tilde{\Gamma} = F_{t^w}(C_{p_d})$ will resemble a fattening of the corresponding graph $\Gamma_Z$. The following is the picture of $\Gamma_Z$ for the standard $Z$.\\
\begin{center}
\begin{picture}(200,180)(-30,0)
\thicklines
\multiput(82,149)(36,0){1}{\line(2,1){18}}
\multiput(82,149)(36,0){1}{\line(-2,1){18}}
\multiput(82,128)(36,0){1}{\line(0,1){21}}
\multiput(46,128)(36,0){2}{\line(2,-1){18}}
\multiput(82,128)(36,0){2}{\line(-2,-1){18}}
\multiput(64,98)(36,0){2}{\line(0,1){21}}
\multiput(28,98)(36,0){3}{\line(2,-1){18}}
\multiput(64,98)(36,0){3}{\line(-2,-1){18}}
\multiput(46,68)(36,0){3}{\line(0,1){21}}
\multiput(10,68)(36,0){4}{\line(2,-1){18}}
\multiput(46,68)(36,0){4}{\line(-2,-1){18}}
\multiput(28,38)(36,0){4}{\line(0,1){21}}
\multiput(-8,38)(36,0){5}{\line(2,-1){18}}
\multiput(28,38)(36,0){5}{\line(-2,-1){18}}
\multiput(10,8)(36,0){5}{\line(0,1){21}}
\put(-26,8){\line(1,0){216}}
\put(-26,8){\line(3,5){108}}
\put(190,8){\line(-3,5){108}}
\end{picture}
\end{center}
\begin{center}
\stepcounter{figure}
Figure \thefigure: standard $\Gamma_Z$ when degree $d=5$
\end{center}
With different convex condition $Z$, the corresponding $\Gamma_Z$ could change to\\
\begin{center}
\begin{picture}(200,200)(-30,0)
\thicklines
\multiput(82,149)(36,0){1}{\line(2,1){18}}
\multiput(82,149)(36,0){1}{\line(-2,1){18}}
\multiput(82,149)(5,-31){2}{\line(2,-3){14}}
\multiput(118,129)(-31,-10){2}{\line(-1,0){23}}
\multiput(86,119)(-27,-10){1}{\line(1,1){10}}
\multiput(46,128)(36,0){1}{\line(2,-1){18}}
\multiput(82,128)(36,0){0}{\line(-2,-1){18}}
\multiput(64,98)(36,0){1}{\line(0,1){21}}
\multiput(28,98)(72,0){2}{\line(2,-1){18}}
\multiput(64,98)(27,-20){2}{\line(1,-2){9}}
\put(73,79){\line(1,0){18}}
\multiput(64,98)(72,0){2}{\line(-2,-1){18}}
\multiput(100,98)(-27,-20){2}{\line(-1,-2){9}}
\multiput(46,68)(72,0){2}{\line(0,1){21}}
\multiput(10,68)(36,0){2}{\line(2,-1){18}}
\multiput(118,68)(36,0){1}{\line(2,-1){18}}
\multiput(46,68)(36,0){1}{\line(-2,-1){18}}
\multiput(118,68)(36,0){2}{\line(-2,-1){18}}
\multiput(28,38)(36,0){4}{\line(0,1){21}}
\multiput(-8,38)(36,0){5}{\line(2,-1){18}}
\multiput(28,38)(36,0){5}{\line(-2,-1){18}}
\multiput(10,8)(36,0){5}{\line(0,1){21}}
\put(-26,8){\line(1,0){216}}
\put(-26,8){\line(3,5){108}}
\put(190,8){\line(-3,5){108}}
\end{picture}
\end{center}
\begin{center}
\stepcounter{figure}
Figure \thefigure: a different $\Gamma_Z$
\end{center}
For details of notations and results, please refer to \cite{N}. In \cite{N}, these results are actually proved in the more general context of curves in arbitrary 2-dimensional toric varieties. After our work was finished, Prof. Oh pointed out to me (during the KIAS conference) the references \cite{M, F} originated from the fundamental work of Gelfand, Kapranov, Zelevinsky \cite{G}, where the authors aimed at very different applications arrived at similar results as in \cite{N} (except our symplectic deformation to graph image). The images of curves under moment map are called amoebas in their work. For more detail on relation to their work, please refer to \cite{N}.\\

The detailed structure of the resulting \l torus fibration of $X_{\psi}$ is described in the following theorem. For detail of the proof, please refer to \cite{lag2,lag3}.\\
\begin{theorem}
\label{ba}
The gradient flow will produce a \l fibration $\tilde{F}: X_{\psi} \rightarrow \partial\Delta$ with singular locus $\tilde{\Gamma} = \tilde{\Gamma}^0 \cup \tilde{\Gamma}^1 \cup \tilde{\Gamma}^2$. There are 4 types of fibres.\\\\
(i). For $p\in \partial\Delta\backslash \tilde{\Gamma}$, $\tilde{F}^{-1}(p)$ is a smooth \l 3-torus.\\
(ii). For $p\in \tilde{\Gamma}^2$, $\tilde{F}^{-1}(p)$ is a \l 3-torus with $2$ circles collapsed to $2$ singular points.\\
(iii). For $p\in \tilde{\Gamma}^1$, $\tilde{F}^{-1}(p)$ is a \l 3-torus with $1$ circles collapsed to $1$ singular points.\\
(iv). For $p\in \tilde{\Gamma}^0$, $\tilde{F}^{-1}(p)$ is a \l 3-torus with $1$ 2-torus collapsed to $1$ singular points.\\
\end{theorem}

In order to get codimension two singular locus, it is necessary to perturb the moment map, which is a Lagrangian fibration, such that the image of the curve under the perturbed moment map is exactly the graph $\Gamma_Z$. This construction is also carried out in \cite{N}. We have\\
\begin{theorem}
\label{da}        
There exists a perturbed Lagrangian fibration $\hat{F}$ of the moment map $F$ satisfying $\hat{F}(C_{s_t}) = \Gamma_Z$.\\
\end{theorem}

With these facts from \cite{N} in hand, general methods developed in \cite{lag1}, \cite{lag2} will enable us to construct Lagrangian torus fibration with one-dimensional singular locus for generic quintic Calabi-Yau near the large complex limit.\\

More precisely, when the generic quintic is near the large complex limit, with the help of theorem \ref{da}, we can produce a Lagrangian fibration $\hat{F}: X_\infty \rightarrow \partial \Delta$ such that $\hat{F}(C)=\Gamma$. $\Gamma$ is a graph in $\partial \Delta \cong S^3$. Its part in each 2-simplex is the kind of graphs described in theorem \ref{da}. Let $\Gamma = \Gamma^1\cup \Gamma^2\cup \Gamma^3$, where $\Gamma^1$ is the smooth part of $\Gamma$, $\Gamma^2$ is the singular part of $\Gamma$ in the interior of the 2-skeleton of $\Delta$, $\Gamma^3$ is the singular part of $\Gamma$ in the 1-skeleton of $\Delta$. Then we have\\
\begin{theorem}
When $X_\psi$ is near the large complex limit, start with Lagrangian fibration $\hat{F}$ the flow of $V$ will produce a Lagrangian fibration $\hat{F}_{\psi}: X_{\psi} \rightarrow \partial\Delta$. There are 4 types of fibres.\\
(i). For $p\in \partial\Delta\backslash \Gamma$, $\hat{F}_{\psi}^{-1}(p)$ is a smooth Lagrangian 3-torus.\\
(ii). For $p\in \Gamma^1$, $\hat{F}_{\psi}^{-1}(p)$ is a type I singular fibre.\\ 
(iii). For $p\in \Gamma^2$, $\hat{F}_{\psi}^{-1}(p)$ is a type II singular fibre.\\
(iv). For $p\in \Gamma^3$, $\hat{F}_{\psi}^{-1}(p)$ is a type III singular fibre.\\
\end{theorem}
Type I singular fibre in the theorem refers to a two-dimensional singular fibre (in this case a nodal ${\bf CP^1}$) times $S^1$. This type of singular fibres have Euler number zero. Other types of singular fibres are illustrated in the following pictures. Notice that type II and III singular fibres have Euler numbers equal to $-1$ and $1$ respectively.\\\\
\begin{center}
\leavevmode
\hbox{%
\epsfxsize=4in
\epsffile{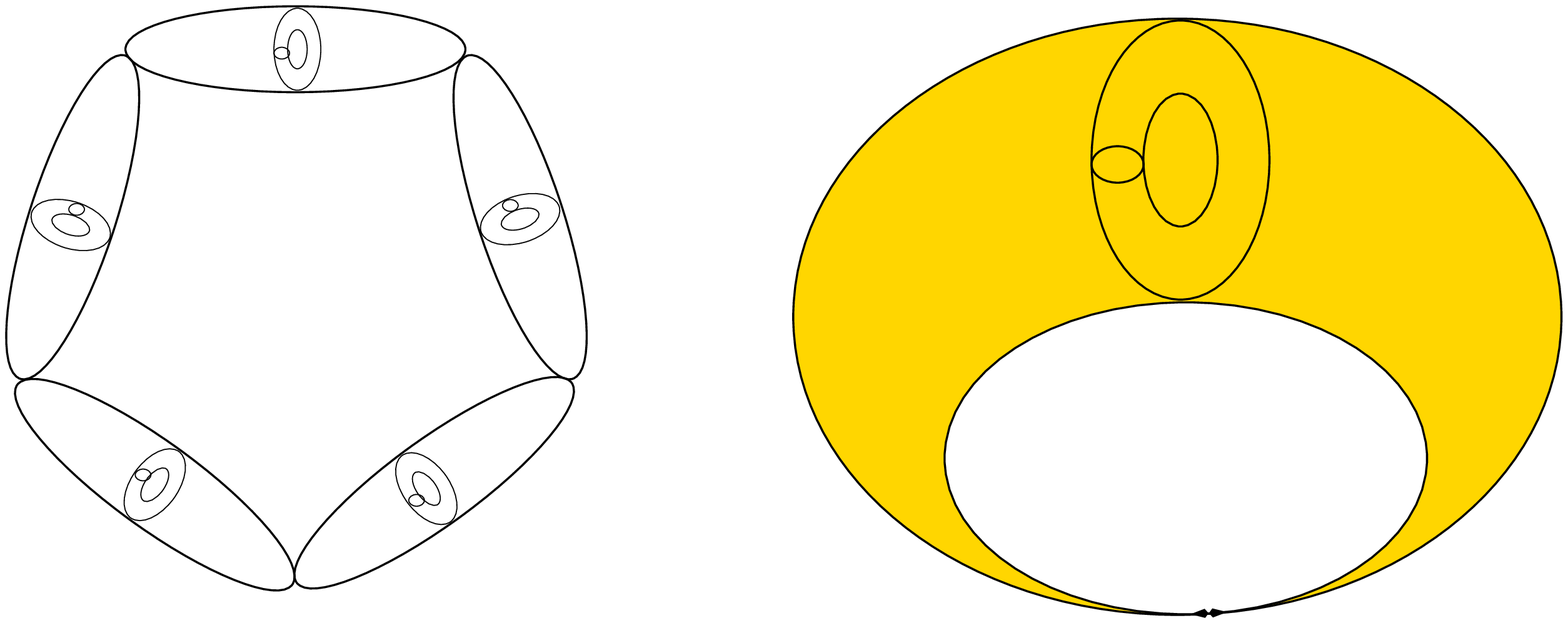}}
\end{center}
\begin{center}
\stepcounter{figure}
Figure \thefigure: Type $III_5$ and type $III$ fibres
\vspace{10pt}
\end{center}
\begin{center}
\leavevmode
\hbox{%
\epsfxsize=4in
\epsffile{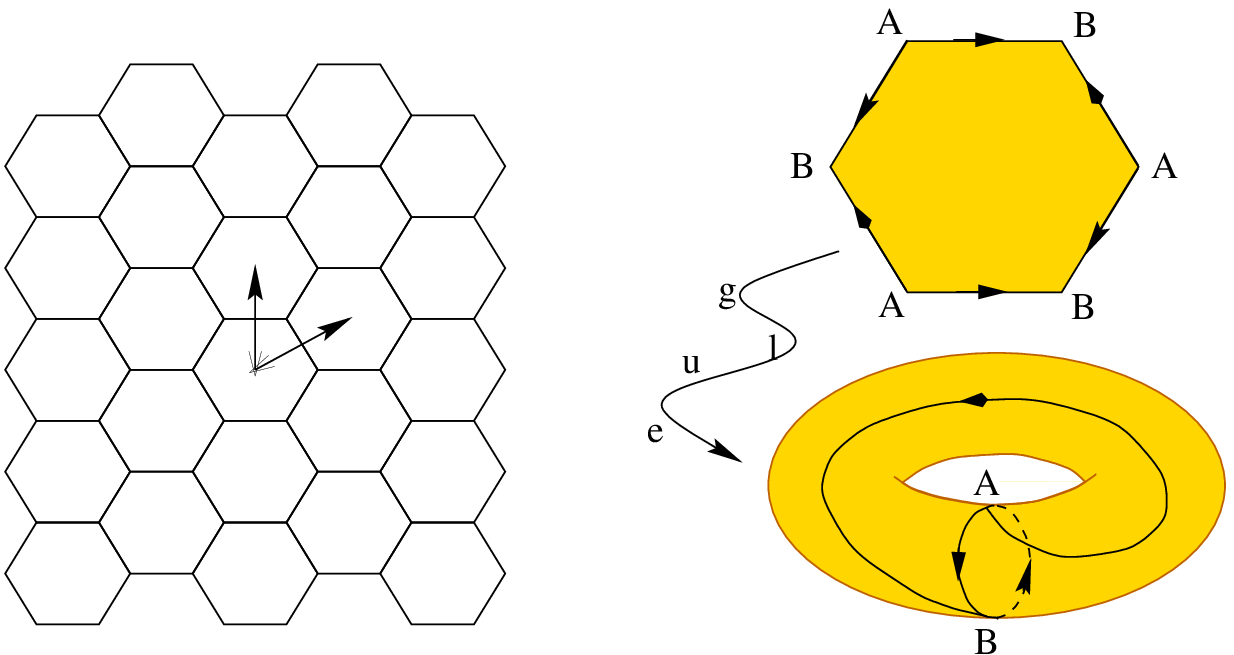}}
\end{center}
\begin{center}
\stepcounter{figure}
Figure \thefigure: Type $II_{5\times 5}$ and type $II$ fibres
\vspace{10pt}
\end{center}
\section{The mirror of quintic}

In the mirror side, the nontrivial part is in the 101 dimensional \k moduli near the large radius limit. As we know the anti-canonical model of the mirror of quintic is the quotient of Fermat type Calabi-Yau\\
\[
Y_\psi = X_\psi/({\bf Z}_5)^3 \subset P_{\Sigma_{\Delta^\vee}} \cong {\bf CP}^4/({\bf Z}_5)^3.
\]\\
Here $\Delta^\vee$ is the dual polyhedron of the polyhedron $\Delta$ corresponding to ${\bf CP}^4$. $\Sigma_{\Delta^\vee}$ is the fan corresponding to $\Delta^\vee$ via normal cone construction. The corresponding toric variety $P_{\Sigma_{\Delta^\vee}}$ is naturally equivalent to ${\bf CP}^4/({\bf Z}_5)^3$.\\

The mirror of quintic \cy are the crepant resolutions of $Y_\psi$. Different crepant resolutions correspond to different chambers of the \k moduli connected by flops. These crepant resolutions of $Y_\psi$ are naturally induced from crepant resolutions of $P_{\Sigma_{\Delta^\vee}} \cong {\bf CP}^4/({\bf Z}_5)^3$, which can be described via crepant subdivision of fan $\Sigma_{\Delta^\vee}$ into different simplical fans. Let $w$ denote the \k class of one such crepant resolution, $w$ determines a simplical fan $\Sigma^w$. $P_{\Sigma^w} \rightarrow P_{\Sigma_{\Delta^\vee}}$ is a crepant resolution. The pullback of $Y_\psi$ (we still use the same notation) is the mirror Calabi-Yau with K\"{a}hler class corresponding to $w$.\\

The gradient flow of the Fermat type quintic Calabi-Yau family $\{X_\psi\}$ is invariant under the action of $({\bf Z}_5)^3$. The quotient gives us the corresponding gradient flow on $P_{\Sigma_{\Delta^\vee}} \cong {\bf CP}^4/({\bf Z}_5)^3$ of the family $\{Y_\psi\}$. This flow pulled back to $P_{\Sigma^w}$ will flow $Y_\infty$ to $Y_\psi$ and induce Lagrangian torus fibration structure on $Y_\psi \subset P_{\Sigma^w}$. Let $F_\psi: Y_\psi \rightarrow \partial \Delta_w$ denote the Lagrangian fibration of $Y_\psi$. Then\\
\begin{prop}
The singular set $C\subset Y_\psi$ of the fibration $F_\psi$ is exactly the intersection of $Y_\psi$ with the complex 2-skeleton of $P_{\Sigma^w}$. In another word,\\
\[
C = Y_\psi \cap {\rm Sing}(Y_\infty).
\]
\end{prop}

The singular locus $\tilde{\Gamma} = F_\psi(C) = F_\infty(C)$ is a fattening of some graph $\Gamma$. We will again use techniques in \cite{lag2, N} to modify $F_\infty(C)$ to $\hat{F}_\infty(C)$ so that $\hat{F}_\infty(C)$ exactly equal to the one-dimensional graph $\Gamma$. To describe this graph $\Gamma$, let's recall from the last section of \cite{lag1} that the singular locus of the Lagrangian fibreation of $Y_\psi \subset P_{\Sigma_{\Delta^\vee}}$ is a fattening of a graph $\hat{\Gamma}\subset \partial \Delta^\vee$, where\\
\[
\hat{\Gamma} = \bigcup_{\{ijklm\}=\{12345\}}\overline{P_{ij}P_{klm}}.
\]
\begin{center}
\leavevmode
\hbox{%
\epsfxsize=4in
\epsffile{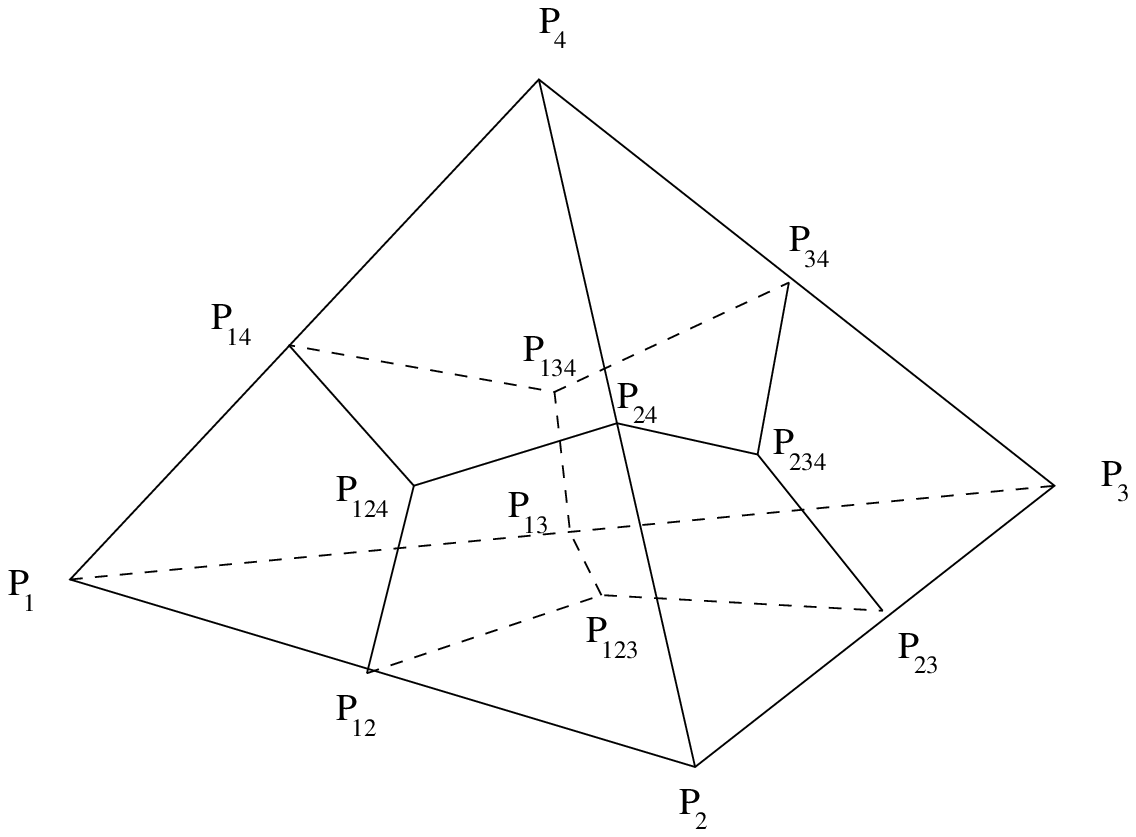}}
\end{center}
\begin{center}
\stepcounter{figure}
Figure \thefigure: $\hat{\Gamma}\subset \partial \Delta^\vee$.
\vspace{10pt}
\end{center}

It is interesting to observe that\\
\[
{\rm Sing}(P_{\Sigma_{\Delta^\vee}}) = {\rm Sing}(Y_\infty).
\]
Hence\\
\[
{\rm Sing}(Y_\psi) = C = Y_\psi \cap {\rm Sing}(Y_\infty).
\]

Let $\tilde{P}_{ij}$ be the point in $C$ that maps to $P_{ij}$ in $\Delta^\vee$. Notice that ${\rm Sing}(C) = \{\tilde{P}_{ij}\}$. Along smooth part of $C$, $Y_\psi$ has $A_5$-singularity. Under the unique crepant resolution, $C$ is turned into 5 copies of $C$. Around $\tilde{P}_{ij} \in {\rm Sing}(C)$, singularity of $Y_\psi$ is much more complicated and crepant resolution is not unique (depending on the \k moduli $w$). The following is a picture (from \cite{lag1}) of fan of such singularity and the subdivision fan of the standard crepant resolution.\\
\begin{center}
\leavevmode
\hbox{%
\epsfxsize=5in
\epsffile{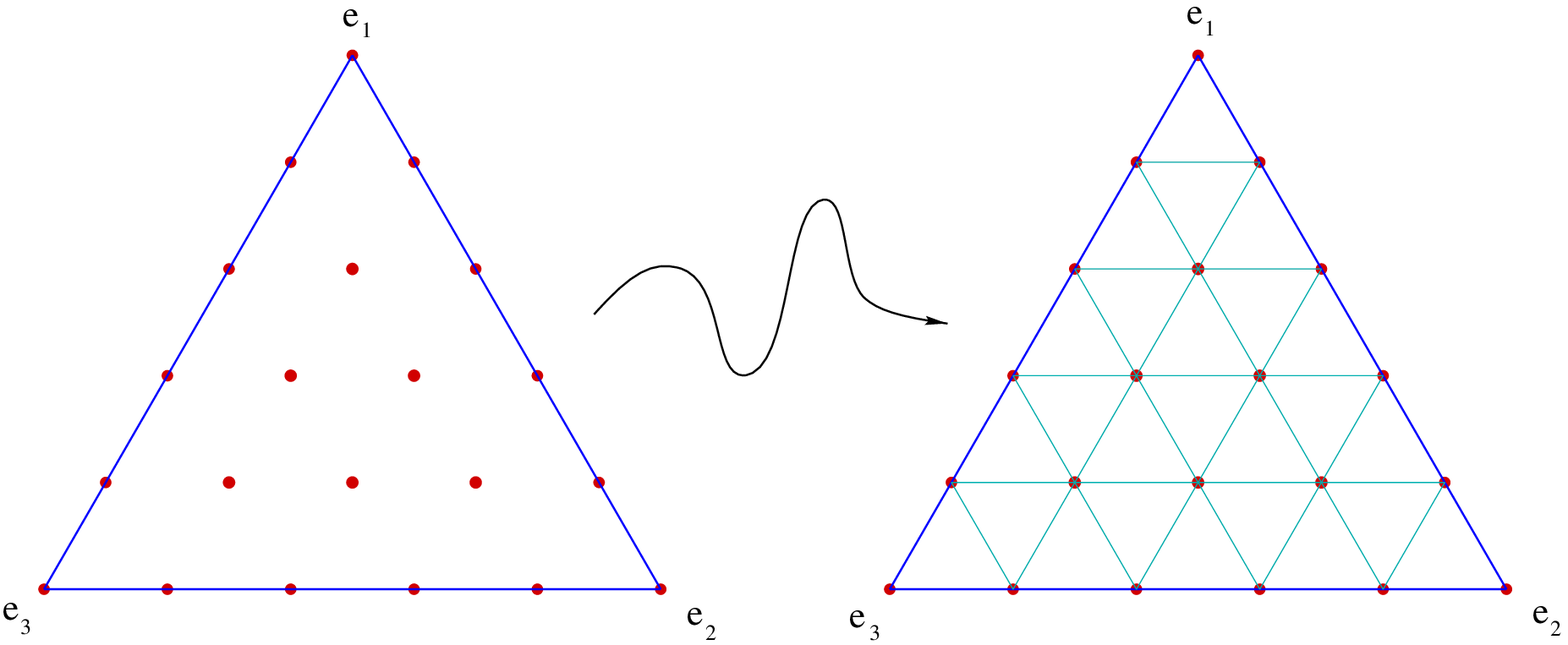}}
\end{center}
\begin{center}
\stepcounter{figure}
Figure \thefigure: the standard crepant resolution of singularity at $\hat{P}_{ij}$
\end{center}

$P_{\Sigma^w} \rightarrow P_{\Sigma_{\Delta^\vee}}$ naturally induces a map $\pi: \Delta_w \rightarrow \Delta^\vee$. We may take $\Gamma$ to be the 1-skeleton of $\pi^{-1}(\hat{\Gamma})$. In a small neighborhood of $P_{ij}$, $\hat{\Gamma} \subset \partial \Delta^\vee$ is indicated in the following picture\\
\begin{center}
\setlength{\unitlength}{1pt}
\begin{picture}(200,100)(-30,10)
\put(80,80){\circle*{2.5}}
\put(90,82){\circle*{2}}
\put(70,78){\circle*{2}}

\multiput(80,80)(36,0){1}{\line(5,1){10}}
\multiput(80,17)(23,0){1}{\line(5,1){10}}
\multiput(26,107)(72,0){1}{\line(5,1){10}}
\multiput(134,107)(36,0){1}{\line(5,1){10}}
\multiput(80,80)(36,0){1}{\line(-5,-1){10}}
\multiput(80,17)(23,0){1}{\line(-5,-1){10}}
\multiput(26,107)(72,0){1}{\line(-5,-1){10}}
\multiput(134,107)(36,0){1}{\line(-5,-1){10}}
\multiput(90,82)(36,0){1}{\line(2,1){54}}
\multiput(90,82)(36,0){1}{\line(-2,1){54}}
\multiput(90,82)(36,0){1}{\line(0,-1){63}}
\multiput(70,78)(36,0){1}{\line(2,1){54}}
\multiput(70,78)(36,0){1}{\line(-2,1){54}}
\multiput(70,78)(36,0){1}{\line(0,-1){63}}
\put(102,48){\vector(-1,0){21}}
\put(110,50){\makebox(0,0){$\hat{\Gamma}$}}
\thicklines
\multiput(80,80)(36,0){1}{\line(2,1){54}}
\multiput(80,80)(36,0){1}{\line(-2,1){54}}
\multiput(80,80)(36,0){1}{\line(0,-1){63}}
\end{picture}
\end{center}
\begin{center}
\stepcounter{figure}
Figure \thefigure: $\hat{\Gamma}\subset \partial \Delta^\vee$ near $P_{ij}$
\end{center}

Under the standard crepant resolution, we get $\Gamma \subset \partial \Delta_w$ as indicated in the following picture\\
\begin{center}
\setlength{\unitlength}{1pt}
\begin{picture}(200,160)(-30,0)
\multiput(78,149)(36,0){1}{\line(5,1){10}}
\multiput(60,119)(36,0){2}{\line(5,1){10}}
\multiput(42,89)(36,0){3}{\line(5,1){10}}
\multiput(24,59)(36,0){4}{\line(5,1){10}}
\multiput(6,29)(36,0){5}{\line(5,1){10}}

\multiput(78,128)(36,0){1}{\line(5,1){10}}
\multiput(60,98)(36,0){2}{\line(5,1){10}}
\multiput(42,68)(36,0){3}{\line(5,1){10}}
\multiput(24,38)(36,0){4}{\line(5,1){10}}

\multiput(6,8)(36,0){5}{\line(5,1){10}}
\multiput(-12,38)(18,30){5}{\line(5,1){10}}
\multiput(96,158)(18,-30){5}{\line(5,1){10}}

\multiput(78,149)(36,0){1}{\line(-5,-1){10}}
\multiput(60,119)(36,0){2}{\line(-5,-1){10}}
\multiput(42,89)(36,0){3}{\line(-5,-1){10}}
\multiput(24,59)(36,0){4}{\line(-5,-1){10}}
\multiput(6,29)(36,0){5}{\line(-5,-1){10}}

\multiput(78,128)(36,0){1}{\line(-5,-1){10}}
\multiput(60,98)(36,0){2}{\line(-5,-1){10}}
\multiput(42,68)(36,0){3}{\line(-5,-1){10}}
\multiput(24,38)(36,0){4}{\line(-5,-1){10}}

\multiput(6,8)(36,0){5}{\line(-5,-1){10}}
\multiput(-12,38)(18,30){5}{\line(-5,-1){10}}
\multiput(96,158)(18,-30){5}{\line(-5,-1){10}}
\thicklines
\multiput(78,149)(36,0){1}{\line(2,1){18}}
\multiput(78,149)(36,0){1}{\line(-2,1){18}}
\multiput(78,128)(36,0){1}{\line(0,1){21}}
\multiput(42,128)(36,0){2}{\line(2,-1){18}}
\multiput(78,128)(36,0){2}{\line(-2,-1){18}}
\multiput(60,98)(36,0){2}{\line(0,1){21}}
\multiput(24,98)(36,0){3}{\line(2,-1){18}}
\multiput(60,98)(36,0){3}{\line(-2,-1){18}}
\multiput(42,68)(36,0){3}{\line(0,1){21}}
\multiput(6,68)(36,0){4}{\line(2,-1){18}}
\multiput(42,68)(36,0){4}{\line(-2,-1){18}}
\multiput(24,38)(36,0){4}{\line(0,1){21}}
\multiput(-12,38)(36,0){5}{\line(2,-1){18}}
\multiput(24,38)(36,0){5}{\line(-2,-1){18}}
\multiput(6,8)(36,0){5}{\line(0,1){21}}

\thinlines
\put(10,2){
\multiput(78,149)(36,0){1}{\line(2,1){18}}
\multiput(78,149)(36,0){1}{\line(-2,1){18}}
\multiput(78,128)(36,0){1}{\line(0,1){21}}
\multiput(42,128)(36,0){2}{\line(2,-1){18}}
\multiput(78,128)(36,0){2}{\line(-2,-1){18}}
\multiput(60,98)(36,0){2}{\line(0,1){21}}
\multiput(24,98)(36,0){3}{\line(2,-1){18}}
\multiput(60,98)(36,0){3}{\line(-2,-1){18}}
\multiput(42,68)(36,0){3}{\line(0,1){21}}
\multiput(6,68)(36,0){4}{\line(2,-1){18}}
\multiput(42,68)(36,0){4}{\line(-2,-1){18}}
\multiput(24,38)(36,0){4}{\line(0,1){21}}
\multiput(-12,38)(36,0){5}{\line(2,-1){18}}
\multiput(24,38)(36,0){5}{\line(-2,-1){18}}
\multiput(6,8)(36,0){5}{\line(0,1){21}}}

\put(-10,-2){
\multiput(78,149)(36,0){1}{\line(2,1){18}}
\multiput(78,149)(36,0){1}{\line(-2,1){18}}
\multiput(78,128)(36,0){1}{\line(0,1){21}}
\multiput(42,128)(36,0){2}{\line(2,-1){18}}
\multiput(78,128)(36,0){2}{\line(-2,-1){18}}
\multiput(60,98)(36,0){2}{\line(0,1){21}}
\multiput(24,98)(36,0){3}{\line(2,-1){18}}
\multiput(60,98)(36,0){3}{\line(-2,-1){18}}
\multiput(42,68)(36,0){3}{\line(0,1){21}}
\multiput(6,68)(36,0){4}{\line(2,-1){18}}
\multiput(42,68)(36,0){4}{\line(-2,-1){18}}
\multiput(24,38)(36,0){4}{\line(0,1){21}}
\multiput(-12,38)(36,0){5}{\line(2,-1){18}}
\multiput(24,38)(36,0){5}{\line(-2,-1){18}}
\multiput(6,8)(36,0){5}{\line(0,1){21}}}
\put(172,18){\vector(-1,0){21}}
\put(179,18){\makebox(0,0){$\Gamma$}}
\end{picture}\\
\end{center}
\begin{center}
\stepcounter{figure}
Figure \thefigure: $\Gamma$ for the standard crepant resolution
\end{center}

For a different crepant resolution, we can get alternative picture for $\Gamma \subset \partial \Delta_w$.\\
\begin{center}
\setlength{\unitlength}{1pt}
\begin{picture}(200,160)(-30,0)
\multiput(78,149)(36,0){1}{\line(5,1){10}}
\multiput(60,119)(23,0){2}{\line(5,1){10}}
\multiput(42,89)(72,0){2}{\line(5,1){10}}
\multiput(24,59)(36,0){4}{\line(5,1){10}}
\multiput(6,29)(36,0){5}{\line(5,1){10}}

\multiput(91,128)(36,0){1}{\line(5,1){10}}
\multiput(60,98)(36,0){2}{\line(5,1){10}}
\multiput(42,68)(72,0){2}{\line(5,1){10}}
\multiput(24,38)(36,0){4}{\line(5,1){10}}
\multiput(69.75,78.5)(16.5,0){2}{\line(5,1){10}}

\multiput(6,8)(36,0){5}{\line(5,1){10}}
\multiput(-12,38)(18,30){5}{\line(5,1){10}}
\multiput(96,158)(18,-30){5}{\line(5,1){10}}

\multiput(78,149)(36,0){1}{\line(-5,-1){10}}
\multiput(60,119)(23,0){2}{\line(-5,-1){10}}
\multiput(42,89)(72,0){2}{\line(-5,-1){10}}
\multiput(24,59)(36,0){4}{\line(-5,-1){10}}
\multiput(6,29)(36,0){5}{\line(-5,-1){10}}

\multiput(91,128)(36,0){1}{\line(-5,-1){10}}
\multiput(60,98)(36,0){2}{\line(-5,-1){10}}
\multiput(42,68)(72,0){2}{\line(-5,-1){10}}
\multiput(24,38)(36,0){4}{\line(-5,-1){10}}
\multiput(69.75,78.5)(16.5,0){2}{\line(-5,-1){10}}

\multiput(6,8)(36,0){5}{\line(-5,-1){10}}
\multiput(-12,38)(18,30){5}{\line(-5,-1){10}}
\multiput(96,158)(18,-30){5}{\line(-5,-1){10}}
\thicklines
\multiput(78,149)(36,0){1}{\line(2,1){18}}
\multiput(78,149)(36,0){1}{\line(-2,1){18}}
\multiput(78,149)(4,-30){2}{\line(2,-3){14}}
\multiput(114,128)(-31,-9){2}{\line(-1,0){23}}
\multiput(82,119)(-27,-10){1}{\line(1,1){10}}
\multiput(42,128)(36,0){1}{\line(2,-1){18}}
\multiput(78,128)(36,0){0}{\line(-2,-1){18}}
\multiput(60,98)(36,0){1}{\line(0,1){21}}
\multiput(24,98)(72,0){2}{\line(2,-1){18}}
\multiput(60,98)(26.5,-19.5){2}{\line(1,-2){9.75}}
\put(69.75,78.5){\line(1,0){16.5}}
\multiput(60,98)(72,0){2}{\line(-2,-1){18}}
\multiput(96,98)(-26.5,-19.5){2}{\line(-1,-2){9.75}}
\multiput(42,68)(72,0){2}{\line(0,1){21}}
\multiput(6,68)(36,0){2}{\line(2,-1){18}}
\multiput(114,68)(36,0){1}{\line(2,-1){18}}
\multiput(42,68)(36,0){1}{\line(-2,-1){18}}
\multiput(114,68)(36,0){2}{\line(-2,-1){18}}
\multiput(24,38)(36,0){4}{\line(0,1){21}}
\multiput(-12,38)(36,0){5}{\line(2,-1){18}}
\multiput(24,38)(36,0){5}{\line(-2,-1){18}}
\multiput(6,8)(36,0){5}{\line(0,1){21}}

\thinlines
\put(10,2){
\multiput(78,149)(36,0){1}{\line(2,1){18}}
\multiput(78,149)(36,0){1}{\line(-2,1){18}}
\multiput(78,149)(4,-30){2}{\line(2,-3){14}}
\multiput(114,128)(-31,-9){2}{\line(-1,0){23}}
\multiput(82,119)(-27,-10){1}{\line(1,1){10}}
\multiput(42,128)(36,0){1}{\line(2,-1){18}}
\multiput(78,128)(36,0){0}{\line(-2,-1){18}}
\multiput(60,98)(36,0){1}{\line(0,1){21}}
\multiput(24,98)(72,0){2}{\line(2,-1){18}}
\multiput(60,98)(26.5,-19.5){2}{\line(1,-2){9.75}}
\put(69.75,78.5){\line(1,0){16.5}}
\multiput(60,98)(72,0){2}{\line(-2,-1){18}}
\multiput(96,98)(-26.5,-19.5){2}{\line(-1,-2){9.75}}
\multiput(42,68)(72,0){2}{\line(0,1){21}}
\multiput(6,68)(36,0){2}{\line(2,-1){18}}
\multiput(114,68)(36,0){1}{\line(2,-1){18}}
\multiput(42,68)(36,0){1}{\line(-2,-1){18}}
\multiput(114,68)(36,0){2}{\line(-2,-1){18}}
\multiput(24,38)(36,0){4}{\line(0,1){21}}
\multiput(-12,38)(36,0){5}{\line(2,-1){18}}
\multiput(24,38)(36,0){5}{\line(-2,-1){18}}
\multiput(6,8)(36,0){5}{\line(0,1){21}}}

\put(-10,-2){
\multiput(78,149)(36,0){1}{\line(2,1){18}}
\multiput(78,149)(36,0){1}{\line(-2,1){18}}
\multiput(78,149)(4,-30){2}{\line(2,-3){14}}
\multiput(114,128)(-31,-9){2}{\line(-1,0){23}}
\multiput(82,119)(-27,-10){1}{\line(1,1){10}}
\multiput(42,128)(36,0){1}{\line(2,-1){18}}
\multiput(78,128)(36,0){0}{\line(-2,-1){18}}
\multiput(60,98)(36,0){1}{\line(0,1){21}}
\multiput(24,98)(72,0){2}{\line(2,-1){18}}
\multiput(60,98)(26.5,-19.5){2}{\line(1,-2){9.75}}
\put(69.75,78.5){\line(1,0){16.5}}
\multiput(60,98)(72,0){2}{\line(-2,-1){18}}
\multiput(96,98)(-26.5,-19.5){2}{\line(-1,-2){9.75}}
\multiput(42,68)(72,0){2}{\line(0,1){21}}
\multiput(6,68)(36,0){2}{\line(2,-1){18}}
\multiput(114,68)(36,0){1}{\line(2,-1){18}}
\multiput(42,68)(36,0){1}{\line(-2,-1){18}}
\multiput(114,68)(36,0){2}{\line(-2,-1){18}}
\multiput(24,38)(36,0){4}{\line(0,1){21}}
\multiput(-12,38)(36,0){5}{\line(2,-1){18}}
\multiput(24,38)(36,0){5}{\line(-2,-1){18}}
\multiput(6,8)(36,0){5}{\line(0,1){21}}}
\put(172,18){\vector(-1,0){21}}
\put(179,18){\makebox(0,0){$\Gamma$}}
\end{picture}
\end{center}
\begin{center}
\stepcounter{figure}
Figure \thefigure: $\Gamma$ for alternative crepant resolution
\end{center}

Notice that such graph changes are caused by flops of the corresponding Calabi-Yau manifolds and locally are of the type\\
\begin{center}
\begin{picture}(200,90)(-150,70)
\thicklines
\put(-50,121){\line(0,1){24}}
\put(-50,121){\line(0,-1){24}}
\put(-50,121){\line(-1,0){30}}
\put(-50,121){\line(1,0){30}}
\put(-50,121){\circle*{4}}
\thinlines
\put(-125,120){\vector(1,0){30}}
\put(-115,0){\thicklines
\put(-59,130){\line(0,1){24}}
\put(-49,120){\line(0,-1){24}}
\put(-59,130){\line(-1,0){30}}
\put(-49,120){\line(1,0){30}}
\put(-59,130){\line(1,-1){10}}
\put(-49,120){\circle*{4}}
\put(-59,130){\circle*{4}}}
\put(-8,121){\vector(1,0){30}}
\put(125,8){\thicklines
\put(-49,122){\line(0,1){24}}
\put(-59,112){\line(0,-1){24}}
\put(-59,112){\line(-1,0){30}}
\put(-49,122){\line(1,0){30}}
\put(-59,112){\line(1,1){10}}
\put(-49,122){\circle*{4}}
\put(-59,112){\circle*{4}}}
\end{picture}
\end{center}
\begin{center}
\stepcounter{figure}
Figure \thefigure: Graph degeneration corresponding to flop
\vspace{10pt}
\end{center}

These singular locus graphs clearly resemble singular locus graphs in figure 2 and 3 obtained via string diagram construction, although the two constructions are quite different.\\

Let $\Gamma = \Gamma^1\cup \Gamma^2\cup \Gamma^3$, where $\Gamma^1$ is the smooth part of $\Gamma$, $\Gamma^3$ is the singular part of $\Gamma$ in the 1-skeleton of $\Delta_w$, $\Gamma^2$ is the rest of singular part of $\Gamma$. With help of some other result in \cite{N}, we can produce a Lagrangian fibration $\hat{F}: Y_\infty \rightarrow \partial \Delta_w$ such that $\hat{F}(C)=\Gamma$. Then we have\\
\begin{theorem}
For $Y_\psi \subset P_{\Sigma^w}$, when $w=(w_m)_{m\in \Delta^0}\in \tau$ is generic and near the large radius limit of $\tau$, start with Lagrangian fibration $\hat{F}$ the flow of $V$ will produce a Lagrangian fibration $\hat{F}_{\psi}: Y_{\psi} \rightarrow \partial\Delta_w$. There are 4 types of fibres.\\
(i). For $p\in \partial\Delta_w\backslash \Gamma$, $\hat{F}_{\psi}^{-1}(p)$ is a smooth Lagrangian 3-torus.\\
(ii). For $p\in \Gamma^1$, $\hat{F}_{\psi}^{-1}(p)$ is a type I singular fibre.\\ 
(iii). For $p\in \Gamma^2$, $\hat{F}_{\psi}^{-1}(p)$ is a type II singular fibre.\\
(iv). For $p\in \Gamma^3$, $\hat{F}_{\psi}^{-1}(p)$ is a type III singular fibre.\\
\end{theorem}

With these constructions of \l torus fibrations for generic quintics and their mirrors, in \cite{lag3}, we were able to prove the symplectic topological version of SYZ conjecture for quintic Calabi-Yau.\\
\begin{theorem}
For generic quintic Calabi-Yau $X$ near the large complex limit, and its mirror Calabi-Yau $Y$ near the large radius limit, there exist corresponding Lagrangian torus fibrations\\
\[
\begin{array}{ccccccc}
X_{s(b)}&\hookrightarrow& X& \ \ &Y_b&\hookrightarrow& Y\\
&&\downarrow& \ &&&\downarrow\\
&& \partial \Delta& \ &&& \partial \Delta_w\\
\end{array}
\]\\
with singular locus $\Gamma \subset \partial \Delta$ and $\Gamma' \subset \partial \Delta_w$, where $s:\partial \Delta_w \rightarrow \partial \Delta$ is a natural homeomorphism with $s(\Gamma')=\Gamma$. For $b\in \partial \Delta_w \backslash \Gamma'$, the corresponding fibres $X_{s(b)}$ and $Y_b$ are naturally dual to each other.\\
\end{theorem}

\section{Calabi-Yau hypersurfaces and complete intersections in toric variety}
In \cite{tor}, we generalize our work in \cite{lag3} to the case of general Calabi-Yau hypersurfaces in toric varieties with respect to reflexive polyhedra, which is exactly the situation of the Batyrev dual polyhedron mirror construction. In a forthcoming paper \cite{ci}, we will further generalize our construction to the case of general Calabi-Yau complete intersections in toric varieties, where mirror construction was proposed by Borisov. Compared to the quintic case, in general toric hypersurface case usually both the \k moduli and the complex moduli of a Calabi-Yau hypersurface are non-trivial. The construction of the Lagrangian torus fibration has to depend on both the \k form and the complex structure of the Calabi-Yau hypersurface. We also need the most general monomial-divisor map to carry out the discussion of symplectic topological SYZ mirror construction for general Calabi-Yau hypersurfaces or complete intersections in toric variety. In the case of complete intersection, The singular locus graph actually exhibit certain knotting phenomenon similar to the flag manifold case we will discuss in the last section. We will only state the result for hypersurfaces here.\\
\begin{theorem}
For generic Calabi-Yau hypersurface $X$ in the toric variety corresponding to a reflexive polyhedron $\Delta$ and its mirror Calabi-Yau hypersurface $Y$ in the toric variety corresponding to the dual reflexive polyhedron $\Delta^\vee$ near their corresponding large complex limit and large radius limit, there exist corresponding Lagrangian torus fibrations\\
\[
\begin{array}{ccccccc}
X_{\phi(b)}&\hookrightarrow& X& \ \ &Y_b&\hookrightarrow& Y\\
&&\downarrow& \ &&&\downarrow\\
&& \partial \Delta_v& \ &&& \partial \Delta_w\\
\end{array}
\]\\
with singular locus $\Gamma \subset \partial \Delta_v$ and $\Gamma' \subset \partial \Delta_w$, where $\phi:\partial \Delta_w \rightarrow \partial \Delta_v$ is a natural homeomorphism, $\phi(\Gamma')=\Gamma$. For $b\in \partial \Delta_w \backslash \Gamma'$, the corresponding fibres $X_{\phi(b)}$ and $Y_b$ are naturally dual to each other. For $b\in \Gamma'$, if $Y_b$ is a singular fibre of type $I$, $II$, $III$, then $X_{\phi(b)}$ is a singular fibre of type $I$, $III$, $II$.\\
\end{theorem}

The original SYZ mirror conjecture was rather sketchy in nature, with no mention of singular locus, singular fibres and duality of singular fibres, which is essential if one wants to use SYZ to construct mirror manifold. Our  discussions on the construction of Lagrangian torus fibrations and symplectic topological SYZ of generic Calabi-Yau hypersurfaces and complete intersections in toric varieties explicitly produce the three types of generic singular fibres (type $I$, $II$, $III$ as described in \cite{lag3}) and exhibit how they are dual to each other under the mirror symmetry. It gives clear indications what should happen in general. In particular, it suggests that type II singular fibre with Euler number $-1$ should be dual to type III singular fibre with Euler number $1$. This together with the knowledge of singular locus from our construction will enable us to give a more precise formulation of SYZ mirror conjecture. This precise formulation  naturally suggests a way to construct mirror manifold from a generic Lagrangian torus fibration of a Calabi-Yau manifold in general. We will state the symplectic version of SYZ conjecture, where duality relation is more precise. In light of examples from Joyce \cite{Joyce}, for the origional special Lagrangian version, singular locus should be a fattening of our graph, singular fibres should also change accordingly, and the identification of singular locus and duality of smooth fibres will be subject to certain quantum effects centered around the large complex (radius) limit, which is not yet fully understood.\\\\
{\bf Precise symplectic SYZ mirror conjecture}\\

For any Calabi-Yau 3-fold $X$, with Calabi-Yau metric $\omega_g$ and holomorphic volume form $\Omega$, there exists a Lagrangian fibration of $X$ over $S^3$\\
\[
\begin{array}{ccc}
T^3&\hookrightarrow& X\\
&&\downarrow\\
&& S^3
\end{array}
\]
with a Lagrangian section and codimension two singular locus $\Gamma\subset S^3$, such that general fibres (over $S^3\backslash \Gamma$) are 3-torus. For generic such fibration, $\Gamma$ is a graph with only 3-valent vertices. Let $\Gamma = \Gamma^1\cup \Gamma^2 \cup \Gamma^3$, where $\Gamma^1$ is the smooth part of $\Gamma$, $\Gamma^2 \cup \Gamma^3$ is the set of the vertices of $\Gamma$. For any leg $\gamma\subset \Gamma^1$, the monodromy of $H_1(X_b)$ of fibre under suitable basis is\\
\[
T_\gamma =\left(
\begin{array}{ccc}
1 &  1 & 0\\
0 &  1 & 0\\
0 &  0 & 1
\end{array}
\right)
\]\\
Singular fibre along $\gamma$ is of type $I$.\\

Consider a vertex $P\in\Gamma^2 \cup \Gamma^3$ with legs $\gamma_1$, $\gamma_2$, $\gamma_3$. Correspondingly, we have monodromy operators $T_1$, $T_2$, $T_3$.\\

For $P\in\Gamma^2$, under suitable basis we have\\
\[
T_1 =\left(
\begin{array}{ccc}
1 &  1 & 0\\
0 &  1 & 0\\
0 &  0 & 1
\end{array}
\right)\ \ 
T_2 =\left(
\begin{array}{ccc}
1 &  0 & -1\\
0 &  1 & 0\\
0 &  0 & 1
\end{array}
\right)\ \ 
T_3 =\left(
\begin{array}{ccc}
1 &  -1 & 1\\
0 &  1 & 0\\
0 &  0 & 1
\end{array}
\right)
\]\\
Singular fibre over $P$ is of type $II$.\\
 
For $P\in\Gamma^3$, under suitable basis we have\\
\[
T_1 =\left(
\begin{array}{ccc}
1 &  0 & 0\\
1 &  1 & 0\\
0 &  0 & 1
\end{array}
\right)\ \ 
T_2 =\left(
\begin{array}{ccc}
1 &  0 & 0\\
0 &  1 & 0\\
-1 &  0 & 1
\end{array}
\right)\ \ 
T_3 =\left(
\begin{array}{ccc}
1 &  0 & 0\\
-1 &  1 & 0\\
1 &  0 & 1
\end{array}
\right)
\]\\
Singular fibre over $P$ is of type $III$.\\

The Lagrangian fibration for the mirror Calabi-Yau manifold $Y$ has the same base $S^3$ and singular locus $\Gamma \subset S^3$. For $b\in S^3\backslash \Gamma$, $Y_b$ is the dual torus of $X_b\cong T^3$. In another word, the $T^3$-fibrations\\
\[
\begin{array}{ccc}
T^3&\hookrightarrow& X\\
&&\downarrow\\
&& S^3\backslash\Gamma
\end{array}
\ \ \ \ \ \ 
\begin{array}{ccc}
T^3&\hookrightarrow& Y\\
&&\downarrow\\
&& S^3\backslash\Gamma
\end{array}
\]\\
are dual to each other. In particular the monodromy operator will be dual to each other.\\

For the fibration of $Y$, singular fibres over $\Gamma^1$ should be type $I$, singular fibres over $\Gamma^2$ should be type $III$, singular fibres over $\Gamma^2$ should be type $II$. Namely, dual singular fibre of a type $I$ singular fibre is still type $I$. Type $II$ and $III$ singular fibres are dual to each other.
\begin{flushright} \rule{2.1mm}{2.1mm} \end{flushright}

Our work in \cite{lag1,lag2,lag3,tor} verified such symplectic SYZ duality for Calabi-Yau hypersurfaces in toric varieties. Our work in progress \cite{ci} will treat the case of complete intersections in toric varieties. \cite{flag} will further push it to the case of Calabi-Yau complete intersections in flag manifolds.\\\\

\section{Conifold transition}

In this section we will discuss the behavior of Lagrangian torus fibrations under conifold transitions of \cy mainifolds. Let $X$ be the \cy manifold that is a smoothing of a \cy conifold $X_0$ with $p$ nodes. Let $\alpha$ be the number of relations among the vanishing cycles in $X$ coming from deformation of the nodes in $X_0$. The conifold transition $Y$ of $X$ is a \cy manifold that is a small resolution of $X_0$. The corresponding birational map $\pi: Y \rightarrow X_0$ contracts $p$ ${\bf CP^1}$'s in $Y$ to $p$ nodes in $X_0$. The topology of $X$ and $Y$ are related as follows \cite{Clemens}.\\
\begin{prop}
\[
h^{1,1}(X) = h^{1,1}(Y) - \alpha.
\]
\[
h^{2,1}(X) = h^{2,1}(Y) + p -\alpha.
\]\\
\end{prop}

We will first describe some local models to indicate the behavior of the Lagrangian fibration near the conifold transition point. Since all these local models have certain $T^2$ or $S^1$ symmetries, it is natural to use symplectic reduction technique to construct fibrations whose fibres are $T^2$ or $S^1$ invariant. Symplectic reduction techniques have been widely used in many fields of mathematics. Recently, Goldstein (\cite{Gold1, Gold2}) applied such techniques to construct certain local special Lagrangian fibrations. Gross in a preprint discussed similar construction. Recently Gross expanded his preprint into \cite{Gross4}, where he also discussed further local examples of generalized special Lagrangian fibrations. For example, among other fibration examples, he briefly discussed the special Lagrangian $T^2\times{\bf R}$ fibration of conifold transition (not as detailed as below) and a diagram similar to Figure 1 of this paper describing corresponding singular locus graph transition. Further discussion of the literature of symplectic reduction techniques can also be found in \cite{Gross4}. (Our gradient flow method also has such flavor.) We are using similar ideas here. Symplectic reduction technique is very powerful when the Calabi-Yau 3-fold has $T^2$ symmetry, while its usefulness will be much limited when the Calabi-Yau 3-fold has merely $S^1$ symmetry. This is probably related to the fact pointed out by Joyce \cite{Joyce} (and also indicated by our gradient flow construction \cite{lag2}) that the special Lagrangian fibrations for such manifolds are not necessarily $C^\infty$.\\\\
{\bf Example:} ($T^2\times{\bf R}$ fibrations) Let\\
\[
X_\epsilon = \{ z=(z_1,z_2,z_3,z_4) \in {\bf C}^4 | p(z) = z_1z_2 - z_3z_4 = \epsilon \}.
\]\\
Then\\
\[
Y = \{ (z, [t_1,t_2]) \in {\bf C}^4\times {\bf CP^1} | t_1z_1 = t_2z_3, t_2z_2 = t_1z_4\}
\]
is a small resolution of $X_0$. Let $\pi: Y \rightarrow X_0$. For $\epsilon > 0$\\
\[
C_1^\epsilon = \{ z=(z_1,z_2,z_3,z_4) \in X_\epsilon | z_2=\bar{z}_1,\ z_4 = -\bar{z}_3\} \cong S^3
\]
is the vanishing cycle.\\
\[
C_2^\epsilon = \{ z=(z_1,z_2,z_3,z_4) \in X_\epsilon | {\rm Re}(z_1)>0,\ {\rm Im}(z_1)=0,\ z_4 = \bar{z}_3\} \cong {\bf R^3}
\]
is the transversal cycle such that $(C_1^\epsilon, C_2^\epsilon)=1$. Let $t = t_2/t_1$, then $\pi^{-1}(C_2^0) \subset Y$ satisfies $x_1\geq 0$ and\\
\[
x_1 = tz_3,\ tx_2 = \bar{z}_3 = z_4.
\]
For fixed $t$, $(x_1,x_2,z_3,\bar{z}_3)$ is uniquely determined up to a positive multiple. Fix $x_1=1$, we have $(1,1/(|t|^2),1/t,1/\bar{t})$. Therefore $\partial (\pi^{-1}(C_2^0)) = \pi^{-1}(0) \cong {\bf CP}^1$.\\

On the other hand, Let\\
\[
S = \{ (z, [1,0]) \in {\bf C}^4\times {\bf CP^1} | z_1 = z_4=0\} \cong {\bf C^2}
\]
be a hypersurface in $Y$ such that $(S, \pi^{-1}(0))=1$. Let $p(z) = f+ ih$ and $\nabla f$ denote the gradient vector field for $f$. Consider\\
\[
V=\frac{\nabla f}{|\nabla f|^2}.
\]
The flow of $V$ will determine a family of symplectic ``blow up'' $\phi_\epsilon: X_\epsilon \rightarrow X_0$. Under the flow, $\pi(S)$ is deformed to symplectic surface $S^\epsilon = \phi_\epsilon^{-1}(\pi(S)) \subset X_\epsilon$ with boundary $\partial (S^\epsilon) = C_1^\epsilon$.\\

Consider a $T^2$ action on $Y$. For $\xi = (\xi_1, \xi_2) \in T^2$, define the action\\
\[
\xi \circ t = \xi_1\xi_2^{-1}t;
\]
\[
\xi \circ z_1 = \xi_1z_1,\ \ \xi \circ z_3 = \xi_2z_3;
\]
\[
\xi \circ z_2 = \xi_1^{-1}z_2,\ \ \xi \circ z_4 = \xi_2^{-1}z_4.
\]
The vector fields corresponding to the two generators are\\
\[
v_1 = 2{\rm Im}\left( z_1\frac{\partial}{\partial z_1} - z_2\frac{\partial}{\partial z_2} + t\frac{\partial}{\partial t}\right),
\]
\[
v_2 = 2{\rm Im}\left( z_3\frac{\partial}{\partial z_3} - z_4\frac{\partial}{\partial z_4} - t\frac{\partial}{\partial t}\right).
\]
The set where the stablizer of $T^2$ is non-trivial is a union of 5 irreducible curves\\
\[
\Delta = \bigcup_{i=0}^4 \Delta_i \subset Y,
\]
where\\
\[
\Delta_0 = \{(z,t)\in Y|z=0\} \cong {\bf CP}^1,
\]
\[
\Delta_i = \{(z,t)\in Y|z_j=0,\ {\rm for}\ j\not=i\} \cong {\bf C}, \ \ {\rm for}\ 1\leq i \leq 4.
\]
Clearly,\\
\[
\Delta_2\cup\Delta_3 = \{(z,t)\in Y|t=0, z_2z_3=0\} \subset \{(z,t)\in Y|t=0\} \cong {\bf C}^2,
\]
\[
\Delta_1\cup\Delta_4 = \{(z,t)\in Y|t=\infty, z_1z_4=0\} \subset \{(z,t)\in Y|t=\infty\} \cong {\bf C}^2.
\]
On $Y$, one can consider the metric\\
\[
\omega_g = \frac{i}{2}\left(\sum_{i=1}^4 dz_i\wedge d\bar{z}_i + \delta \frac{dt\wedge d\bar{t}}{(1+|t|^2)^2}\right).
\]
We can compute
\[
\i(v_1) \omega_g = \frac{1}{2}d\left(|z_1|^2 - |z_2|^2 - \frac{\delta}{1+|t|^2}\right),
\]
\[
\i(v_2) \omega_g = \frac{1}{2}d\left(|z_3|^2 - |z_4|^2 + \frac{\delta}{1+|t|^2}\right).
\]
$Y/T^2$ can be naturally identified with ${\bf R}^4$ via\\
\[
\left(|z_1|^2-|z_2|^2 - \frac{\delta}{1+|t|^2}, |z_3|^2-|z_4|^2 + \frac{\delta}{1+|t|^2}, z_1z_2+z_3z_4\right).
\]\\
$\rho: Y \rightarrow {\bf R}^3$ defined by\\
\[
\rho(z) = \left(|z_1|^2-|z_2|^2 - \frac{\delta}{1+|t|^2}, |z_3|^2-|z_4|^2 + \frac{\delta}{1+|t|^2}, {\rm Re}(z_1z_2+z_3z_4)\right)
\]
is a $T^2\times {\bf R}$ fibration. The singular set of the fibration is exactly $\Delta$. The singular locus $\rho(\Delta) \subset \{\rho_3 =0\} \cong {\bf R}^2$ is a graph as follows.\\
\begin{center}
\begin{picture}(200,170)(-150,30)
\thicklines
\put(-59,130){\line(0,1){72}}
\put(-49,166){\makebox(0,0){$\Gamma_3$}}
\put(-41,112){\line(0,-1){72}}
\put(-31,76){\makebox(0,0){$\Gamma_4$}}
\put(-59,130){\line(-1,0){108}}
\put(-113,140){\makebox(0,0){$\Gamma_2$}}
\put(-41,112){\line(1,0){108}}
\put(13,122){\makebox(0,0){$\Gamma_1$}}
\put(-59,130){\line(1,-1){18}}
\put(-44,127){\makebox(0,0){$\Gamma_0$}}
\put(-41,112){\circle*{4}}
\put(-59,130){\circle*{4}}
\end{picture}
\end{center}
\begin{center}
\stepcounter{figure}
Figure \thefigure: Singular locus of the $T^2\times {\bf R}$ fibration of $Y$
\vspace{10pt}
\end{center}
When $\delta$ shrinks to 0, $Y$ is blown down to $X_0$ and the singular locus of the fibration changes to\\
\begin{center}
\begin{picture}(200,170)(-150,30)
\thicklines
\put(-50,121){\line(0,1){72}}
\put(-40,157){\makebox(0,0){$\Gamma_3$}}
\put(-50,121){\line(0,-1){72}}
\put(-40,85){\makebox(0,0){$\Gamma_4$}}
\put(-50,121){\line(-1,0){108}}
\put(-104,131){\makebox(0,0){$\Gamma_2$}}
\put(-50,121){\line(1,0){108}}
\put(4,131){\makebox(0,0){$\Gamma_1$}}
\put(-50,121){\circle*{4}}
\end{picture}
\end{center}
\begin{center}
\stepcounter{figure}
Figure \thefigure: Singular locus of the $T^2\times {\bf R}$ fibration of $X_0$

\vspace{10pt}
\end{center}
The $T^2$ action on $Y$ can naturally be carried over to a $T^2$ action on $X_\epsilon$, defined as\\
\[
\xi \circ z_1 = \xi_1z_1,\ \ \xi \circ z_3 = \xi_2z_3;
\]
\[
\xi \circ z_2 = \xi_1^{-1}z_2,\ \ \xi \circ z_4 = \xi_2^{-1}z_4.
\]
$X_\epsilon/T^2$ can be naturally identified with ${\bf R}^4$ via $(|z_1|^2-|z_2|^2, |z_3|^2-|z_4|^2, z_1z_2+z_3z_4)$.\\
\[
\rho_\epsilon: X_\epsilon \rightarrow {\bf R}^3
\]
defined by\\
\[
\rho(z) = (|z_1|^2-|z_2|^2, |z_3|^2-|z_4|^2, {\rm Re}(z_1z_2+z_3z_4))
\]
is a $T^2\times {\bf R}$ fibration. The singular set is now union of 2 curves\\
\[
\Delta^\epsilon = \Delta^\epsilon_{12} \cup \Delta^\epsilon_{34}
\]
where\\
\[
\Delta^\epsilon_{12} = \{z\in X_\epsilon|z_3=z_4=0,\ z_1z_2=\epsilon\} \cong {\bf C}^*,
\]
\[
\Delta^\epsilon_{34} = \{z\in X_\epsilon|z_1=z_2=0,\ z_3z_4=-\epsilon\} \cong {\bf C}^*.
\]
The singular locus $\rho(\Delta^\epsilon)$ consists of two lines not in a plane.\\
\begin{center}
\begin{picture}(200,200)(-150,30)
\put(-59,130){\line(1,-1){18}}
\put(-59,211){\line(1,-1){18}}
\put(-59,67){\line(1,-1){18}}
\put(-59,67){\line(0,1){144}}
\multiput(-131,49)(216,0){2}{\line(0,1){144}}
\multiput(-131,49)(0,144){2}{\line(1,0){216}}
\multiput(-149,67)(216,0){2}{\line(0,1){144}}
\multiput(-149,67)(0,144){2}{\line(1,0){216}}
\thicklines
\put(-41,112){\line(0,1){16}}
\put(-41,132){\line(0,1){61}}
\put(-41,112){\line(0,-1){63}}
\put(-31,100){\makebox(0,0){$\Gamma_{34}$}}
\put(-59,130){\line(-1,0){90}}
\put(-59,130){\line(1,0){126}}
\put(-4,140){\makebox(0,0){$\Gamma_{12}$}}
\end{picture}
\end{center}
\begin{center}
\stepcounter{figure}
Figure \thefigure: Singular locus of the $T^2\times {\bf R}$ fibration of $X_\epsilon$
\vspace{10pt}
\end{center}
\begin{flushright} \rule{2.1mm}{2.1mm} \end{flushright}
It is not hard to verify that this fibration is actually the so-called generalized special Lagrangian fibration.\\\\
{\bf Example:} ($S^1\times{\bf R}^2$ fibrations)\\
Consider the same situation as in the previous example\\
\[
X_\epsilon = \{ z=(z_1,z_2,z_3,z_4) \in {\bf C}^4 | p(z) = z_1z_2 - z_3z_4 = \epsilon \}
\]
\[
Y = \{ (z, [t_1,t_2]) \in {\bf C}^4\times {\bf CP^1} | t_1z_1 = t_2z_3, t_2z_2 = t_1z_4\}
\]
is a small resolution of $X_0$.\\\\
Consider an $S^1$ action on $Y$. For $\xi \in S^1$, define the action\\
\[
\xi \circ t = \xi^{-1} t;
\]
\[
\xi \circ z_1 = z_1,\ \ \xi \circ z_3 = \xi z_3;
\]
\[
\xi \circ z_2 = z_2,\ \ \xi \circ z_4 = \xi^{-1} z_4.
\]
The vector field corresponding to the generator is\\
\[
v = 2{\rm Im}\left( z_3\frac{\partial}{\partial z_3} - z_4\frac{\partial}{\partial z_4} - t\frac{\partial}{\partial t}\right).
\]
The set where the stablizer of $S^1$ is non-trivial is a union of 2 irreducible curves\\
\[
\Delta = \Delta_1 \cup \Delta_2 \subset Y,
\]
where\\
\[
\Delta_i = \{(z,t)\in Y|z_j=0,\ {\rm for}\ j\not=i\} \cong {\bf C}, \ \ {\rm for}\ 1\leq i \leq 2.
\]
On $Y$, one can consider the metric\\
\[
\omega_g = \frac{i}{2}\left(\sum_{i=1}^4 dz_i\wedge d\bar{z}_i + \delta \frac{dt\wedge d\bar{t}}{(1+|t|^2)^2}\right).
\]
We have
\[
\i(v) \omega_g = \frac{1}{2}d\left(|z_3|^2 - |z_4|^2 + \frac{\delta}{1+|t|^2}\right).
\]
$Y/S^1$ can be naturally identified with ${\bf R}^5$ via\\
\[
\left(z_1, z_2, |z_3|^2-|z_4|^2 + \frac{\delta}{1+|t|^2}\right).
\]\\
$\rho: Y \rightarrow {\bf R}^3$ defined by\\
\[
\rho(z) = \left({\rm Re}(z_1), {\rm Re}(z_2), |z_3|^2-|z_4|^2 + \frac{\delta}{1+|t|^2}\right)
\]
is an $S^1\times {\bf R}^2$ fibration. The singular set of the fibration is exactly $\Delta$. The singular locus $\Gamma = \rho(\Delta) = \Gamma_1 \cup \Gamma_2 \subset {\bf R}^2$ consists of two lines not in a plane ($\Gamma_1 = \rho(\Delta_1) = \{\rho_2 =\rho_3 =0\}$ and $\Gamma_2 = \rho(\Delta_2) = \{\rho_1 =\rho_3 =0\}$).\\
\begin{center}
\begin{picture}(200,200)(-150,30)
\put(-59,130){\line(1,-1){18}}
\put(-59,211){\line(1,-1){18}}
\put(-59,67){\line(1,-1){18}}
\put(-59,67){\line(0,1){144}}
\multiput(-131,49)(216,0){2}{\line(0,1){144}}
\multiput(-131,49)(0,144){2}{\line(1,0){216}}
\multiput(-149,67)(216,0){2}{\line(0,1){144}}
\multiput(-149,67)(0,144){2}{\line(1,0){216}}
\thicklines
\put(-41,112){\line(0,1){16}}
\put(-41,132){\line(0,1){61}}
\put(-41,112){\line(0,-1){63}}
\put(-31,100){\makebox(0,0){$\Gamma_2$}}
\put(-59,130){\line(-1,0){90}}
\put(-59,130){\line(1,0){126}}
\put(-4,140){\makebox(0,0){$\Gamma_1$}}
\end{picture}
\end{center}
\begin{center}
\stepcounter{figure}
Figure \thefigure: Singular locus of the $S^1\times {\bf R}^2$ fibration of $Y$
\vspace{10pt}
\end{center}
When $\delta$ shrinks to 0, $Y$ is blown down to $X_0$ and the singular locus of the fibration changes to\\
\begin{center}
\begin{picture}(200,170)(-150,30)
\thicklines
\put(-50,121){\vector(0,1){72}}
\put(-40,157){\makebox(0,0){$\Gamma_2$}}
\put(-40,195){\makebox(0,0){$\rho_2$}}
\put(-50,121){\line(0,-1){72}}
\put(-40,85){\makebox(0,0){$\Gamma_2$}}
\put(-50,121){\vector(-1,0){108}}
\put(-104,131){\makebox(0,0){$\Gamma_1$}}
\put(-161,131){\makebox(0,0){$\rho_1$}}
\put(-50,121){\line(1,0){108}}
\put(4,131){\makebox(0,0){$\Gamma_1$}}
\put(-50,121){\circle*{4}}
\end{picture}
\end{center}
\begin{center}
\stepcounter{figure}
Figure \thefigure: Singular locus of the $S^1\times {\bf R}^2$ fibration of $X_0$
\vspace{10pt}
\end{center}
The $S^1$ action on $Y$ can naturally be carried over to an $S^1$ action on $X_\epsilon$, defined as\\
\[
\xi \circ z_1 = z_1,\ \ \xi \circ z_3 = \xi z_3;
\]
\[
\xi \circ z_2 = z_2,\ \ \xi \circ z_4 = \xi^{-1}z_4.
\]
$X_\epsilon/S^1$ can be naturally identified with ${\bf R}^5$ via $(z_1, z_2, |z_3|^2-|z_4|^2)$.\\
\[
\rho_\epsilon: X_\epsilon \rightarrow {\bf R}^3
\]
defined by\\
\[
\rho(z) = ({\rm Re}(z_1), {\rm Re}(z_2), |z_3|^2-|z_4|^2)
\]
is an $S^1\times {\bf R}^2$ fibration. The singular set is now\\
\[
\Delta^\epsilon = \{z\in X_\epsilon|z_3=z_4=0,\ z_1z_2=\epsilon\} \cong {\bf C}^*.
\]\\
Assume that $\epsilon$ is positive, then the singular locus of $\rho$ is\\
\[
\tilde{\Gamma} = \rho(\Delta^\epsilon) = \{0 \leq \rho_1\rho_2 \leq \epsilon,\rho_3=0,\} \subset \{\rho_3=0\}.
\]
It contains two parts\\
\[
\tilde{\Gamma}_+ = \{\rho_1 \geq 0,\rho_2 \geq 0,\rho_1\rho_2 \leq \epsilon,\rho_3=0\},
\]
\[
\tilde{\Gamma}_- = \{\rho_1 \leq 0,\rho_2 \leq 0,\rho_1\rho_2 \leq \epsilon,\rho_3=0\}.
\]
The vanishing circle of $\Delta^\epsilon$ is\\
\[
S^1_\epsilon = \{(\sqrt{\epsilon}e^{i\theta}, \sqrt{\epsilon}e^{i\theta}, 0, 0)\}.
\]
Under the fibration map\\
\[
\Gamma_0 = \rho(S^1_\epsilon) = \{(t, t, 0)| -1\leq t \leq 1\}.
\]
The following is a picture of $\tilde{\Gamma}$.\\
\begin{center}
\setlength{\unitlength}{1pt}
\begin{picture}(300,300)(-150,-150)
\multiput(-18,0)(-18,0){7}{\vector(0,1){10}}
\multiput(18,0)(18,0){7}{\vector(0,-1){10}}
\multiput(0,-18)(0,-18){7}{\vector(1,0){10}}
\multiput(0,18)(0,18){7}{\vector(-1,0){10}}
\put(-36,36){\vector(1,-1){10}}
\put(-54,24){\vector(1,-1){10}}
\put(-72,18){\vector(1,-1){10}}
\put(-90,15){\vector(1,-1){10}}
\put(-108,12){\vector(1,-1){10}}
\put(-126,12){\vector(1,-1){10}}
\put(-24,54){\vector(1,-1){10}}
\put(-18,72){\vector(1,-1){10}}
\put(-15,90){\vector(1,-1){10}}
\put(-12,108){\vector(1,-1){10}}
\put(-12,126){\vector(1,-1){10}}

\put(36,-36){\vector(-1,1){10}}
\put(54,-24){\vector(-1,1){10}}
\put(72,-18){\vector(-1,1){10}}
\put(90,-15){\vector(-1,1){10}}
\put(108,-12){\vector(-1,1){10}}
\put(126,-12){\vector(-1,1){10}}
\put(24,-54){\vector(-1,1){10}}
\put(18,-72){\vector(-1,1){10}}
\put(15,-90){\vector(-1,1){10}}
\put(12,-108){\vector(-1,1){10}}
\put(12,-126){\vector(-1,1){10}}
\thicklines
\put(-42,30){\line(1,1){12}}
\put(-30,42){\line(1,2){9}}
\put(-21,60){\line(1,5){9}}
\put(-12,105){\line(0,1){25}}
\put(-42,30){\line(-2,-1){18}}
\put(-60,21){\line(-5,-1){45}}
\put(-105,12){\line(-1,0){25}}

\put(42,-30){\line(-1,-1){12}}
\put(30,-42){\line(-1,-2){9}}
\put(21,-60){\line(-1,-5){9}}
\put(12,-105){\line(0,-1){25}}
\put(42,-30){\line(2,1){18}}
\put(60,-21){\line(5,1){45}}
\put(105,-12){\line(1,0){25}}

\put(0,0){\vector(0,1){144}}
\put(10,146){\makebox(0,0){$\rho_2$}}
\put(0,0){\line(0,-1){144}}
\put(0,0){\vector(-1,0){180}}
\put(-188,2){\makebox(0,0){$\rho_1$}}
\put(0,0){\line(1,0){180}}
\put(0,0){\circle*{4}}
\put(-18,18){\makebox(0,0){${\bf \tilde{\Gamma}_+}$}}
\put(18,-18){\makebox(0,0){${\bf \tilde{\Gamma}_-}$}}
\end{picture}
\end{center}
\begin{center}
\stepcounter{figure}
Figure \thefigure: Singular locus of the $S^1\times {\bf R}^2$ fibration of $X_\epsilon$
\vspace{10pt}
\end{center}
Deforming the fibration along the arrows, we can get an $S^1\times {\bf R}^2$ fibration with singular locus\\
\begin{center}
\begin{picture}(200,160)(-150,30)
\thicklines
\put(-59,130){\line(0,1){72}}
\put(-49,166){\makebox(0,0){$\Gamma_+$}}
\put(-41,112){\line(0,-1){72}}
\put(-31,76){\makebox(0,0){$\Gamma_-$}}
\put(-59,130){\line(-1,0){108}}
\put(-113,140){\makebox(0,0){$\Gamma_+$}}
\put(-41,112){\line(1,0){108}}
\put(13,122){\makebox(0,0){$\Gamma_-$}}
\put(-59,130){\line(1,-1){18}}
\put(-44,127){\makebox(0,0){$\Gamma_0$}}
\put(-41,112){\circle*{4}}
\put(-59,130){\circle*{4}}
\end{picture}
\end{center}
\begin{center}
\stepcounter{figure}
Figure \thefigure: Modified singular locus of the $S^1\times {\bf R}^2$ fibration of $X_\epsilon$
\vspace{10pt}
\end{center}
\begin{flushright} \rule{2.1mm}{2.1mm} \end{flushright}

The singular set of the fibration in a fibre over $\Gamma_+$ and $\Gamma_-$ is a line. The singular set of the fibration in a fibre over $\Gamma_0$ is a union of two parallel lines. The singular set of the fibration in a fibre over the 3-valent points of $\Gamma$ is a graph identical to the graph $\Gamma$.\\

Unfortunately, these $S^1\times {\bf R}^2$ fibrations are not even Lagrangian fibrations. In \cite{con}, we will discuss this local conifold transition for more general metrics and modify these explicit non-Lagrangian fibrations into explicit Lagrangian ones with the same topological structure.\\

From these $S^1\times {\bf R}^2$ fibrations, it is also easy to see why we required the components of the singular locus to be located in parallel planes as mentioned in the introduction. Over the singular locus, the invariant $S^1$ is vanishing, which implies that the components of the singular locus will be located in some level planes of $\rho_3 = |z_3|^2 -|z_4|^2$.\\\\

\section{Calabi-Yau complete intersection in flag manifold}

In \cite{BCKS1,BCKS2}, Batyrev {\it et. al.} were able to generalize the mirror construction to the case of Calabi-Yau complete intersections in partial flag manifolds, which had much to do with conifolds transitions. It is interesting to see how we can generalize our construction of Lagrangian fibration to this case and what new phenomenon might ocurr.\\

With the above local models in mind, using the gradient flow approach together with some symplectic patching technique, we will be able to construct Lagrangian torus fibration for Calabi-Yau complete intersections in flag manifolds and their mirror, therefore proving the symplectic SYZ mirror symmetry for such classes of Calabi-Yau. In this section we will illustrate the key ideas through an example (the case of Grassmannian ${\bf Gr(2,4)}$). As one might notice, ${\bf Gr(2,4)}$ is actually a conic hypersurface in ${\bf CP}^5$. The Calabi-Yau hypersurface in ${\bf Gr(2,4)}$ is actually a complete intersection in toric variety. In the Lagrangian torus fibration point of view, this is a rare coincidence. As we will show in \cite{ci,flag}, Lagrangian torus fibration structure for Calabi-Yau complete intersections in toric varieties and Calabi-Yau hypersurfaces in flag manifolds have very distinct characters. In fact, for this example, the large complex limits from the flag manifold point of view and from the toric point of view are entirely different in conventional sense. Consequently the Lagrangian torus fibrations constructed from the two points of view are also entirely different. In any case, our method here will not really rely on the fact that ${\bf Gr(2,4)}$ is a conic hypersurface in ${\bf CP}^5$, and key ideas for the general case already show up here. We will provide details of the general case in the forthcoming paper \cite{flag}.\\\\
{\bf Example:} Consider the case of ${\bf Gr(2,4)}$, $z \in \Lambda^2{\bf H_C}$ can be expressed as\\
\[
z = \sum_{1\leq i\leq j \leq 4} z_{ij} e_i\wedge e_j.
\]
\[
z\wedge z = (z_{12}z_{34} - z_{13}z_{24} + z_{14}z_{23})e_1\wedge e_2 \wedge e_3 \wedge e_4.
\]
${\bf Gr(2,4)} \subset {\bf P(\Lambda^2 H_C)} \cong {\bf CP}^5$ is a quadric defined by\\
\[
z_{12}z_{34} - z_{13}z_{24} + z_{14}z_{23} =0.
\]\\
Let $P(2,4)\subset {\bf CP}^5$ be the 4-dimensional Gorenstein toric Fano variety defined by the quadratic equation\\
\[
z_{13}z_{24} - z_{14}z_{23} = 0.
\]\\
Then $P(2,4)$ is a degeneration of ${\bf Gr(2,4)}$ in terms of the family\\
\[
z_{13}z_{24} - z_{14}z_{23} - \epsilon z_{12}z_{34} = 0.
\]\\

Consider the fan description of ${\bf CP}^5$ with homogeneous coordinates $[z_{jk}]_{1\leq j < k \leq 4}$. Then $M$ is naturally\\
\[
M = \left\{z^I=\prod_{1\leq j < k \leq 4}z_{jk}^{i_{jk}}\ \left|\ I=(i_{jk})_{1\leq j < k \leq 4} \in{\bf Z}^6,\ \ |I| = \sum_{1\leq j < k \leq 4}i_{jk} = 0\right.\right\}.
\]
Let
\[
\Delta = \left\{z^I \in M \ \left|\ I + I_0 \geq 0\right.\right\}
\]\\
where $I_0 = (1,1,\cdots,1)$. Let\\
\[
s_{m_1} = z_{13}z_{24},\ \ s_{m_2} = z_{14}z_{23}.
\]\\
Then $P(2,4) \cong \{s = s_{m_1} - s_{m_2} =0\}$. $m_1-m_2 = (i_{jk})_{1\leq j < k \leq 4} \in{\bf Z}^6$, where $i_{12}=i_{34}=0$, $i_{13}=i_{24}=1$ and $i_{14}=i_{23}=-1$. Naturally\\
\[
N = \left.\left\{w=(w_{jk})_{1\leq j < k \leq 4} \in{\bf Z}^6\right\}\right/ {\bf Z}\cdot w_0,
\]\\
where $w_0 = (1,1,\cdots,1)$.\\
\[
\Delta^\vee = \left\{[w] \in N \ \left|\ w \geq 0,\ |w|\leq 1\right.\right\}.
\]
\[
M_s = M/(m_1-m_2).
\]
\[
N_s = (m_1-m_2)^\perp = \left\{[w] \in N \ \left|\ w_{13}+w_{24}=w_{14}+w_{23}\right.\right\}.
\]\\
1-cones in $\Sigma_s$ come from 1-cones and 2-cones in $\Sigma$ that intersect $(m_1-m_2)^\perp$. Assume $w = (w_{12},w_{34},w_{13},w_{24},w_{14},w_{23})$. The generating vectors are\\
\[
w_1 = (1,0,0,0,0,0),\ w_3 = (0,0,1,0,1,0),\ w_5 = (0,0,0,1,1,0), 
\]
\[
w_2 = (0,1,0,0,0,0),\ w_4 = (0,0,1,0,0,1),\ w_6 = (0,0,0,1,0,1). 
\]\\
These are the vertices of $\Delta_s^\vee$.\\

Up to symmetries, there are 3 distinguished 1-simplices $\overline{w_1w_2}$, $\overline{w_1w_3}$, $\overline{w_3w_4}$. There are no additional integral points in any one of them. Up to symmetries, there are 3 distinguished 2-faces $\overline{w_1w_2w_3}$, $\overline{w_1w_3w_4}$, $\overline{w_3w_4w_6w_5}$. There are no additional integral points in any one of them. Up to symmetries, there are 2 distinguished 3-faces $\overline{w_1w_2w_3w_4}$, $\overline{w_1w_3w_4w_6w_5}$. There are no additional integral points in any one of them.\\

The fan $\Sigma_s$ corresponding to $\Delta_s \subset M_s$ consists of cones over subfaces of $\Delta_s^\vee$. Recall that\\
\[
P(2,4) \cong P_{\Sigma_s} = \bigcup_{\sigma \in \Sigma_s} T_\sigma.
\]\\
The singular set of $P_{\Sigma_s}$ is a line $l = \{z_{11} = z_{13} = z_{14} = z_{23} =0\}\cong{\bf CP}^1$ which is an union of $T_\sigma$ corresponding to non-simplicial cones $\sigma$ over $\overline{w_1w_3w_4w_6w_5}$, $\overline{w_3w_4w_6w_5}$, $\overline{w_2w_3w_4w_6w_5}$. Singularity along $l$ is 3-dimensional conifold singularity characterized by the cone over $\overline{w_3w_4w_6w_5}$.\\

Assume $I = (i_{12},i_{34},i_{13},i_{24},i_{14},i_{23})$. The vertices of $\Delta_s$ are\\
\[
I_1 = (1,-3,1,1,0,0),\ I_3 = (1,1,1,-3,0,0),\ I_5 = (1,1,0,0,1,-3), 
\]
\[
I_2 = (-3,1,1,1,0,0),\ I_4 = (1,1,0,0,-3,1),\ I_6 = (1,1,-3,1,0,0). 
\]\\
Under the moment map, $l$ is mapped to $\overline{I_1I_2}$. Small resolution of $l$ will result in a modified polyhedron $\hat{\Delta}_s$ of $\Delta_s$, where edge $\overline{I_1I_2}$ is replaced by a 2-dimensional parallelgram.\\

A Calabi-Yau hypersurface $X$ in $P(2,4)$ is the intersection of a quartic and $P(2,4)$ in ${\bf CP^5}$. Using the gradient flow approach developed in \cite{lag1,lag2,lag3,tor}, when $X$ is near the large complex limit, we can construct a Lagrangian torus fibration of $X$ over $\hat{\Delta}_s$ with codimension 2 singular locus $\Gamma$. We will describe $\Delta_s$ and $\Gamma$ in detail.\\
 
Up to symmetries, there are 3 distinguished 1-simplices $\overline{I_1I_2}$, $\overline{I_1I_3}$, $\overline{I_3I_4}$. On $\overline{I_1I_2}$, there are 3 other integral points\\
\[
(0,-2,1,1,0,0),\ (-1,-1,1,1,0,0),\ (-2,0,1,1,0,0).
\] 
On $\overline{I_1I_3}$, there are 3 other integral points\\
\[
(1,-2,1,0,0,0),\ (1,-1,1,-1,0,0),\ (1,0,1,-2,0,0).
\]
On $\overline{I_3I_4}$, there are 3 other integral points\\
\[
(1,1,1,-2,-1,0),\ (1,1,1,-1,-2,0),\ (1,1,0,-1,-2,1).
\]

Up to symmetries, there are 3 distinguished 2-faces $\overline{I_1I_2I_3}$, $\overline{I_1I_3I_4}$, $\overline{I_3I_4I_6I_5}$. On $\overline{I_1I_2I_3}$ (or $\overline{I_1I_3I_4}$), there are 15 integral points. Graph $\Gamma$ in $\overline{I_1I_2I_3}$ is determined by certain simplicial decomposition of $\overline{I_1I_2I_3}$. For the standard symplicial decomposition, we have the following graph $\Gamma$\\\\
\begin{center}
\setlength{\unitlength}{1.1pt}
\begin{picture}(200,180)(0,45)
\put(118,226){\makebox(0,0){$I_1$}}
\put(212,62){\makebox(0,0){$I_3$}}
\put(-4,62){\makebox(0,0){$I_2$}}
\put(80,184){\line(1,0){48}}
\put(152,64){\line(3,5){24}}
\put(56,64){\line(-3,5){24}}
\put(56,144){\line(1,0){96}}
\put(104,64){\line(3,5){48}}
\put(104,64){\line(-3,5){48}}
\put(32,104){\line(1,0){144}}
\put(56,64){\line(3,5){72}}
\put(152,64){\line(-3,5){72}}
\put(8,64){\line(1,0){192}}
\put(8,64){\line(3,5){96}}
\put(200,64){\line(-3,5){96}}
\multiput(104,224)(48,0){1}{\circle*{4}}
\multiput(80,184)(48,0){2}{\circle*{4}}
\multiput(56,144)(48,0){3}{\circle*{4}}
\multiput(32,104)(48,0){4}{\circle*{4}}
\multiput(8,64)(48,0){5}{\circle*{4}}
\thicklines
\multiput(104,198)(48,0){1}{\line(2,1){24}}
\multiput(104,198)(48,0){1}{\line(-2,1){24}}
\multiput(104,170)(48,0){1}{\line(0,1){28}}
\multiput(56,170)(48,0){2}{\line(2,-1){24}}
\multiput(104,170)(48,0){2}{\line(-2,-1){24}}
\multiput(80,130)(48,0){2}{\line(0,1){28}}
\multiput(32,130)(48,0){3}{\line(2,-1){24}}
\multiput(80,130)(48,0){3}{\line(-2,-1){24}}
\multiput(56,90)(48,0){3}{\line(0,1){28}}
\multiput(8,90)(48,0){4}{\line(2,-1){24}}
\multiput(56,90)(48,0){4}{\line(-2,-1){24}}
\multiput(32,50)(48,0){4}{\line(0,1){28}}
\end{picture}
\end{center}
\begin{center}
\stepcounter{figure}
Figure \thefigure: Possible singular locus in $\overline{I_1I_2I_3}$
\vspace{10pt}
\end{center}
If the simplicial decomposition is changed, $\Gamma$ will change accordingly\\
\begin{center}
\setlength{\unitlength}{1.1pt}
\begin{picture}(200,180)(0,45)
\put(118,226){\makebox(0,0){$I_1$}}
\put(212,62){\makebox(0,0){$I_3$}}
\put(-4,62){\makebox(0,0){$I_2$}}
\put(82,184){\line(5,-3){70}}
%\put(82,181){\line(2,-1){68}}
\put(80,184){\line(1,0){48}}
\put(152,64){\line(3,5){24}}
\put(56,64){\line(-3,5){24}}
\put(56,144){\line(1,0){96}}
\put(104,64){\line(3,5){48}}
\put(104,64){\line(-3,5){48}}
\put(32,104){\line(1,0){48}}
\put(128,104){\line(1,0){48}}
\put(104,64){\line(0,1){80}}
\put(56,64){\line(3,5){48}}
\put(152,64){\line(-3,5){72}}
\put(8,64){\line(1,0){192}}
\put(8,64){\line(3,5){96}}
\put(200,64){\line(-3,5){96}}
\multiput(104,224)(48,0){1}{\circle*{4}}
\multiput(80,184)(48,0){2}{\circle*{4}}
\multiput(56,144)(48,0){3}{\circle*{4}}
\multiput(32,104)(48,0){4}{\circle*{4}}
\multiput(8,64)(48,0){5}{\circle*{4}}
\thicklines
\put(104,198){\line(2,1){24}}
\put(104,198){\line(-2,1){24}}
\put(104,198){\line(2,-3){18}}
\put(109,158){\line(2,-3){19}}
\put(152,171){\line(-1,0){31}}
\put(80,158){\line(1,0){30}}
\put(56,170){\line(2,-1){24}}
\put(109,158){\line(1,1){13}}
\put(80,130){\line(0,1){28}}
\multiput(32,130)(96,0){2}{\line(2,-1){24}}
\multiput(80,130)(35,-26){2}{\line(1,-2){13}}
\put(93,104){\line(1,0){22}}
\multiput(80,130)(96,0){2}{\line(-2,-1){24}}
\multiput(128,130)(-35,-26){2}{\line(-1,-2){13}}
\multiput(56,90)(96,0){2}{\line(0,1){28}}
\multiput(8,90)(48,0){2}{\line(2,-1){24}}
\put(152,90){\line(2,-1){24}}
\put(56,90){\line(-2,-1){24}}
\multiput(152,90)(48,0){2}{\line(-2,-1){24}}
\multiput(32,50)(48,0){4}{\line(0,1){28}}
\end{picture}
\end{center}
\begin{center}
\stepcounter{figure}
Figure \thefigure: Possible singular locus in $\overline{I_1I_2I_3}$
\vspace{10pt}
\end{center}
On $\overline{I_3I_4I_6I_5}$, there are 25 integral points. For the standard symplicial decomposition, we have the following graph $\Gamma$\\ 
\begin{center}
\setlength{\unitlength}{1.4pt}
\begin{picture}(200,160)(-30,-10)
\put(-10,144){\makebox(0,0){$I_3$}}
\put(-10,0){\makebox(0,0){$I_4$}}
\put(154,144){\makebox(0,0){$I_5$}}
\put(154,0){\makebox(0,0){$I_6$}}
\multiput(0,144)(36,0){5}{\line(0,-1){144}}
\multiput(0,144)(0,-36){5}{\line(1,0){144}}
\multiput(0,108)(36,0){4}{\line(1,1){36}}
\multiput(0,72)(36,0){4}{\line(1,1){36}}
\multiput(0,36)(36,0){4}{\line(1,1){36}}
\multiput(0,0)(36,0){4}{\line(1,1){36}}
\multiput(0,144)(36,0){5}{\circle*{4}}
\multiput(0,108)(36,0){5}{\circle*{4}}
\multiput(0,72)(36,0){5}{\circle*{4}}
\multiput(0,36)(36,0){5}{\circle*{4}}
\multiput(0,0)(36,0){5}{\circle*{4}}
\thicklines
\multiput(12,132)(36,0){5}{\line(-2,-1){24}}
\multiput(12,132)(36,0){4}{\line(1,2){12}}
\multiput(12,132)(36,0){4}{\line(1,-1){12}}
\multiput(12,96)(36,0){5}{\line(-2,-1){24}}
\multiput(12,96)(36,0){4}{\line(1,2){12}}
\multiput(12,96)(36,0){4}{\line(1,-1){12}}
\multiput(12,60)(36,0){5}{\line(-2,-1){24}}
\multiput(12,60)(36,0){4}{\line(1,2){12}}
\multiput(12,60)(36,0){4}{\line(1,-1){12}}
\multiput(12,24)(36,0){5}{\line(-2,-1){24}}
\multiput(12,24)(36,0){4}{\line(1,2){12}}
\multiput(12,24)(36,0){4}{\line(1,-1){12}}
\multiput(12,-12)(36,0){4}{\line(1,2){12}}
\end{picture}
\end{center}
\begin{center}
\stepcounter{figure}
Figure \thefigure: Possible singular locus in $\overline{I_3I_4I_6I_5}$
\vspace{10pt}
\end{center}
Of course, suitable changes of simplicial decomposition will result in changes of corresponding $\Gamma$.\\

Singular set $l$ of $P_{\Sigma_s}$ is a ${\bf CP}^1$ corresponding to $\overline{I_1I_2}$. Around $\overline{I_1I_2}$, before the small resolution, $\Gamma$ looks like\\
\begin{center}
\begin{picture}(200,160)(50,-10)
\put(-10,0){\makebox(0,0){$I_1$}}
\put(298,144){\makebox(0,0){$I_2$}}
\put(0,0){\line(2,1){288}}
\multiput(0,0)(72,36){5}{\circle*{4}}
\thicklines
\multiput(36,18)(72,36){4}{\line(0,1){28}}
\multiput(36,18)(72,36){4}{\line(0,-1){28}}
\multiput(36,18)(72,36){4}{\line(1,0){28}}
\multiput(36,18)(72,36){4}{\line(-1,0){28}}
\end{picture}
\end{center}
\begin{center}
\stepcounter{figure}
Figure \thefigure: Singular locus around $\overline{I_1I_2}$
\vspace{10pt}
\end{center}
$\Gamma$ in $\overline{I_1I_2I_3}$, $\overline{I_1I_2I_4}$, $\overline{I_1I_2I_5}$, $\overline{I_1I_2I_6}$ intersect at 4 points in $\overline{I_1I_2}$. When one resolve singular set $l$ by small resolution, then $\Gamma$ around $\overline{I_1I_2}$ changes to\\
\begin{center}
\begin{picture}(200,160)(50,-10)
\put(-8,-6){\makebox(0,0){$I_1$}}
\put(296,150){\makebox(0,0){$I_2$}}
\put(5,-5){\line(2,1){288}}
\put(-5,5){\line(2,1){288}}
\multiput(-5,5)(72,36){5}{\line(1,-1){10}}
\multiput(-5,5)(72,36){5}{\circle*{4}}
\multiput(5,-5)(72,36){5}{\circle*{4}}
\thicklines
\multiput(31,23)(72,36){4}{\line(0,1){28}}
\multiput(41,13)(72,36){4}{\line(0,-1){28}}
\multiput(41,13)(72,36){4}{\line(1,0){28}}
\multiput(31,23)(72,36){4}{\line(-1,0){28}}
\multiput(31,23)(72,36){4}{\line(1,-1){10}}
\end{picture}
\end{center}
\begin{center}
\stepcounter{figure}
Figure \thefigure: Singular locus around $\overline{I_1I_2}$ after small resolution
\vspace{10pt}
\end{center}
Instead, if one smooths out the singular set $l$ of $P_{\Sigma_s} \subset {\bf CP}^5$ by deformation in ${\bf CP}^5$, then $\Gamma$ around $\overline{I_1I_2}$ changes to\\ 
\begin{center}
\begin{picture}(200,160)(50,-10)
\put(-10,0){\makebox(0,0){$I_1$}}
\put(298,144){\makebox(0,0){$I_2$}}
\put(0,0){\line(2,1){288}}
\multiput(31,23)(72,36){4}{\line(1,-1){10}}
\multiput(0,0)(72,36){5}{\circle*{4}}
\thicklines
\multiput(41,24)(72,36){4}{\line(0,1){22}}
\multiput(41,22)(72,36){4}{\line(0,-1){30}}
\multiput(31,23)(72,36){4}{\line(1,0){32}}
\multiput(31,23)(72,36){4}{\line(-1,0){20}}
\end{picture}
\end{center}
\begin{center}
\stepcounter{figure}
Figure \thefigure: Singular locus around $\overline{I_1I_2}$ after smoothing
\vspace{10pt}
\end{center}
$\Gamma$ in $\overline{I_1I_2I_3} \cup \overline{I_1I_2I_4} \cup \overline{I_1I_2I_5} \cup \overline{I_1I_2I_6}$, which used to intersect at 4 points in $\overline{I_1I_2}$, now form two pieces $\overline{I_1I_6I_2I_3}$ and $\overline{I_1I_5I_2I_4}$ (like in the following picture) that avoid each other and intertwine near $\overline{I_1I_2}$ as in the above picture.\\\\
\begin{center}
\setlength{\unitlength}{1.1pt}
\begin{picture}(200,200)(-40,40)
\put(118,226){\makebox(0,0){$I_1$}}
\put(212,62){\makebox(0,0){$I_3$}}
\put(-4,62){\makebox(0,0){$I_2$}}
\put(-100,222){\makebox(0,0){$I_6$}}
\put(80,184){\line(1,0){48}}
\put(152,64){\line(3,5){24}}
\put(56,64){\line(-3,5){24}}
\put(56,144){\line(1,0){96}}
\put(104,64){\line(3,5){48}}
\put(104,64){\line(-3,5){48}}
\put(32,104){\line(1,0){144}}
\put(56,64){\line(3,5){72}}
\put(152,64){\line(-3,5){72}}
\put(8,64){\line(1,0){192}}
\put(8,64){\line(3,5){96}}
\put(200,64){\line(-3,5){96}}
\multiput(104,224)(48,0){1}{\circle*{4}}
\multiput(80,184)(48,0){2}{\circle*{4}}
\multiput(56,144)(48,0){3}{\circle*{4}}
\multiput(32,104)(48,0){4}{\circle*{4}}
\multiput(8,64)(48,0){5}{\circle*{4}}

\put(-16,104){\line(1,0){48}}
\put(56,224){\line(3,-5){24}}
\put(-40,224){\line(-3,-5){24}}
\put(-40,144){\line(1,0){96}}
\put(8,224){\line(3,-5){48}}
\put(8,224){\line(-3,-5){48}}
\put(-64,184){\line(1,0){144}}
\put(-40,224){\line(3,-5){72}}
\put(56,224){\line(-3,-5){72}}
\put(-88,224){\line(1,0){192}}
\put(-88,224){\line(3,-5){96}}
\put(104,224){\line(-3,-5){96}}
\multiput(8,64)(48,0){1}{\circle*{4}}
\multiput(-16,104)(48,0){2}{\circle*{4}}
\multiput(-40,144)(48,0){3}{\circle*{4}}
\multiput(-64,184)(48,0){4}{\circle*{4}}
\multiput(-88,224)(48,0){5}{\circle*{4}}

\thicklines
\multiput(104,198)(48,0){1}{\line(2,1){24}}
\multiput(104,198)(48,0){1}{\line(-2,1){24}}
\multiput(104,170)(48,0){1}{\line(0,1){28}}
\multiput(56,170)(48,0){2}{\line(2,-1){24}}
\multiput(104,170)(48,0){2}{\line(-2,-1){24}}
\multiput(80,130)(48,0){2}{\line(0,1){28}}
\multiput(32,130)(48,0){3}{\line(2,-1){24}}
\multiput(80,130)(48,0){3}{\line(-2,-1){24}}
\multiput(56,90)(48,0){3}{\line(0,1){28}}
\multiput(8,90)(48,0){4}{\line(2,-1){24}}
\multiput(56,90)(48,0){4}{\line(-2,-1){24}}
\multiput(32,50)(48,0){4}{\line(0,1){28}}

\multiput(8,90)(48,0){1}{\line(2,-1){24}}
\multiput(8,90)(48,0){1}{\line(-2,-1){24}}
\multiput(8,118)(48,0){1}{\line(0,-1){28}}
\multiput(-40,118)(48,0){2}{\line(2,1){24}}
\multiput(8,118)(48,0){2}{\line(-2,1){24}}
\multiput(-16,158)(48,0){2}{\line(0,-1){28}}
\multiput(-64,158)(48,0){3}{\line(2,1){24}}
\multiput(-16,158)(48,0){3}{\line(-2,1){24}}
\multiput(-40,198)(48,0){3}{\line(0,-1){28}}
\multiput(-88,198)(48,0){4}{\line(2,1){24}}
\multiput(-40,198)(48,0){4}{\line(-2,1){24}}
\multiput(-64,238)(48,0){4}{\line(0,-1){28}}

\end{picture}
\end{center}
\begin{center}
\stepcounter{figure}
Figure \thefigure: Singular locus in $\overline{I_1I_6I_2I_3}$
\vspace{10pt}
\end{center}
\begin{center}
\setlength{\unitlength}{1.1pt}
\begin{picture}(200,200)(-40,40)
\put(118,226){\makebox(0,0){$I_1$}}
\put(212,62){\makebox(0,0){$I_4$}}
\put(-4,62){\makebox(0,0){$I_2$}}
\put(-100,222){\makebox(0,0){$I_5$}}
\put(82,184){\line(5,-3){70}}
%\put(82,181){\line(2,-1){68}}
\put(80,184){\line(1,0){48}}
\put(152,64){\line(3,5){24}}
\put(56,64){\line(-3,5){24}}
\put(56,144){\line(1,0){96}}
\put(104,64){\line(3,5){48}}
\put(104,64){\line(-3,5){48}}
\put(32,104){\line(1,0){48}}
\put(128,104){\line(1,0){48}}
\put(104,64){\line(0,1){80}}
\put(56,64){\line(3,5){48}}
\put(152,64){\line(-3,5){72}}
\put(8,64){\line(1,0){192}}
\put(8,64){\line(3,5){96}}
\put(200,64){\line(-3,5){96}}
\multiput(104,224)(48,0){1}{\circle*{4}}
\multiput(80,184)(48,0){2}{\circle*{4}}
\multiput(56,144)(48,0){3}{\circle*{4}}
\multiput(32,104)(48,0){4}{\circle*{4}}
\multiput(8,64)(48,0){5}{\circle*{4}}

\put(-14,104){\line(5,3){70}}
\put(-16,104){\line(1,0){48}}
\put(56,224){\line(3,-5){24}}
\put(-40,224){\line(-3,-5){24}}
\put(-40,144){\line(1,0){96}}
\put(8,224){\line(3,-5){48}}
\put(8,224){\line(-3,-5){48}}
\put(-64,184){\line(1,0){48}}
\put(32,184){\line(1,0){48}}
\put(8,224){\line(0,-1){80}}
\put(-40,224){\line(3,-5){48}}
\put(56,224){\line(-3,-5){72}}
\put(-88,224){\line(1,0){192}}
\put(-88,224){\line(3,-5){96}}
\put(104,224){\line(-3,-5){96}}
\multiput(8,64)(48,0){1}{\circle*{4}}
\multiput(-16,104)(48,0){2}{\circle*{4}}
\multiput(-40,144)(48,0){3}{\circle*{4}}
\multiput(-64,184)(48,0){4}{\circle*{4}}
\multiput(-88,224)(48,0){5}{\circle*{4}}
\thicklines
\put(104,198){\line(2,1){24}}
\put(104,198){\line(-2,1){24}}
\put(104,198){\line(2,-3){18}}
\put(109,158){\line(2,-3){19}}
\put(152,171){\line(-1,0){31}}
\put(80,158){\line(1,0){30}}
\put(56,170){\line(2,-1){24}}
\put(109,158){\line(1,1){13}}
\put(80,130){\line(0,1){28}}
%\multiput(32,130)(96,0){1}{\line(2,-1){24}}
\multiput(128,130)(96,0){1}{\line(2,-1){24}}
\multiput(80,130)(35,-26){2}{\line(1,-2){13}}
\put(93,104){\line(1,0){22}}
\multiput(80,130)(96,0){2}{\line(-2,-1){24}}
\multiput(128,130)(-35,-26){2}{\line(-1,-2){13}}
\multiput(56,90)(96,0){2}{\line(0,1){28}}
\multiput(8,90)(48,0){2}{\line(2,-1){24}}
\put(152,90){\line(2,-1){24}}
\put(56,90){\line(-2,-1){24}}
\multiput(152,90)(48,0){2}{\line(-2,-1){24}}
\multiput(32,50)(48,0){4}{\line(0,1){28}}

\put(8,90){\line(2,-1){24}}
\put(8,90){\line(-2,-1){24}}
\put(8,90){\line(2,3){18}}
\put(13,130){\line(2,3){19}}
\put(56,117){\line(-1,0){31}}
\put(-16,130){\line(1,0){30}}
\put(-40,118){\line(2,1){24}}
\put(13,130){\line(1,-1){13}}
\put(-16,158){\line(0,-1){28}}
\multiput(-64,158)(96,0){2}{\line(2,1){24}}
\multiput(-16,158)(35,26){2}{\line(1,2){13}}
\put(-3,184){\line(1,0){22}}
\multiput(-16,158)(96,0){2}{\line(-2,1){24}}
\multiput(32,158)(-35,26){2}{\line(-1,2){13}}
\multiput(-40,198)(96,0){2}{\line(0,-1){28}}
\multiput(-88,198)(48,0){2}{\line(2,1){24}}
\put(56,198){\line(2,1){24}}
\put(-40,198){\line(-2,1){24}}
\multiput(56,198)(48,0){2}{\line(-2,1){24}}
\multiput(-64,238)(48,0){4}{\line(0,-1){28}}
\end{picture}
\end{center}
\begin{center}
\stepcounter{figure}
Figure \thefigure: Singular locus in $\overline{I_1I_5I_2I_4}$

\vspace{10pt}
\end{center}

One can easily observe similarity of these two graphs with the graph for $\overline{I_3I_4I_6I_5}$. Totally, $\Gamma$ consists of 3 square pieces $\Gamma_{3465}$, $\Gamma_{1623}$, $\Gamma_{1524}$ and 8 triangle pieces $\Gamma_{134}$, $\Gamma_{145}$, $\Gamma_{156}$, $\Gamma_{163}$, $\Gamma_{234}$, $\Gamma_{245}$, $\Gamma_{256}$, $\Gamma_{263}$. This will give us the singular locus of the Lagrangian torus fibrations of Calabi-Yau hypersurfaces in the Grassmannian ${\bf Gr(2,4)}$. Notice that unlike the hypersurfaces in toric varieties, where the singular locus graph is completely unknotted, in the Grassmannian (and more generally the flag manifold) case, the singular locus graphs show slight knotting phenomenon. In fact such knotting phenomenon already show up in the case of complete intersections in toric varieties (in somewhat different fashion), which we will discuss in our forthcoming paper \cite{ci}. It is a rather interesting question to see if more complicated knotting is allowed for the singular locus graph of Lagrangian torus fibrations of more general Calabi-Yau manifolds.\\

In summary, the Calabi-Yau hypersurface in $\hat{P}(2,4)$ is transformed to Calabi-Yau hypersurface in the Grassmannian ${\bf Gr(2,4)}$ through conifold transition at 4 ordinary double points in the singular Calabi-Yau hypersurface in $P(2,4)$. The corresponding singular locus graph locally undergoes the following basic graph change as mentioned in the introduction.\\
\begin{center}
\begin{picture}(200,90)(-150,70)
\thicklines
\put(-50,121){\line(0,1){24}}
\put(-50,121){\line(0,-1){24}}
\put(-50,121){\line(-1,0){30}}
\put(-50,121){\line(1,0){30}}
\put(-50,121){\circle*{4}}
\thinlines
\put(-125,120){\vector(1,0){30}}
\put(-115,0){\thicklines
\put(-59,130){\line(0,1){24}}
\put(-49,120){\line(0,-1){24}}
\put(-59,130){\line(-1,0){30}}
\put(-49,120){\line(1,0){30}}
\put(-59,130){\line(1,-1){10}}
\put(-49,120){\circle*{4}}
\put(-59,130){\circle*{4}}}
\put(-8,121){\vector(1,0){30}}
\put(125,-5){
\put(-57,128){\circle*{1}}
\put(-55,126){\circle*{1}}
\put(-53,124){\circle*{1}}
\put(-51,122){\circle*{1}}
\thicklines
\put(-49,120){\line(0,1){8}}
\put(-49,132){\line(0,1){20}}
\put(-49,120){\line(0,-1){20}}
\put(-59,130){\line(-1,0){30}}
\put(-59,130){\line(1,0){42}}}
\end{picture}
\end{center}
\begin{center}
\stepcounter{figure}
Figure \thefigure: Graph degeneration corresponding to conifold transition
\vspace{10pt}
\end{center}
Singular locus graph and the singular fibres are modelled according to the local $T^2\times {\bf R}$ fibration examples discussed in the previous section.\\

We can similarly discuss the Lagrangian torus fibrations in the mirror side. The mirror of the Calabi-Yau hypersurface in $\hat{P}(2,4)$ is transformed to the mirror of Calabi-Yau hypersurface in the Grassmannian ${\bf Gr(2,4)}$ also through conifold transition at 4 ordinary double points (of course in reverse direction). Singular locus graph and the singular fibres are modelled according to the local $S^1\times {\bf R}^2$ fibration examples discussed in the previous section. Symplectic SYZ duality can also be proved in this situation. We hope to provide more details of these arguments in \cite{con,flag}.\\\\
{\bf Acknowledgement:} I would like to thank Prof. S.T. Yau for constant encouragement, Prof. Kefeng Liu for suggesting the work of \cite{BCKS1,BCKS2} to me during the conference in Montreal, Prof. Yong-Geun Oh for pointing out the work of \cite{M} to me during the KIAS conference and Prof. A. Libgober for helpful discussion. I would also like to thank Prof. Yong-Geun Oh for inviting me to the well organized conference on mirror symmetry held in KIAS. Major part of this work was done while I was in Columbia University and finished up in University of Illinois at Chicago. I am very grateful to both universities for excellent research environment. Thanks also go to Qin Jing for stimulating discussions and suggestions.\\\\

\end{document}